\documentclass[a4paper,oneside,english]{elsart}
\usepackage{ae}
\usepackage{aecompl}
\usepackage[T1]{fontenc}
\usepackage[latin1]{inputenc}
\usepackage{array}
\usepackage{longtable}
\usepackage{subfigure}
\usepackage{amsmath}
\usepackage{graphicx}
\usepackage{amssymb}
\IfFileExists{url.sty}{\usepackage{url}}
                      {\newcommand{\url}{\texttt}}

\makeatletter

\providecommand{\tabularnewline}{\\}

 \newenvironment{lyxlist}[1]
   {\begin{list}{}
     {\settowidth{\labelwidth}{#1}
      \setlength{\leftmargin}{\labelwidth}
      \addtolength{\leftmargin}{\labelsep}
      }}
   {\end{list}}

\usepackage{subfigure}
\numberwithin{equation}{section}
\date{22nd December 2005}
\usepackage{babel}
\makeatother
\begin{document}
\begin{frontmatter}

\title{Iterative algorithms to approximate canonical Gabor windows: Computational
aspects}

\author{A.J.E.M. Janssen}

\address{Philips Research Laboratories WO-02, 5656AA Eindhoven, The Netherlands}

\ead{a.j.e.m.janssen@philips.com}

\author{Peter L. Søndergaard\corauthref{cor}}

\corauth[cor]{Corresponding author.}

\address{Technical University of Denmark, Department of Mathematics, Building
303, 2800 Lyngby, Denmark.}

\ead{P.Soendergaard@mat.dtu.dk}

\begin{abstract}
In this paper we investigate the computational aspects of some recently
proposed iterative methods for approximating the canonical tight and
canonical dual window of a Gabor frame $\left(g,a,b\right)$. The
iterations start with the window $g$ while the iteration steps comprise
the window $g$, the $k^{th}$ iterand $\gamma_{k}$, the frame operators
$S$ and $S_{k}$ corresponding to $\left(g,a,b\right)$ and $\left(\gamma_{k},a,b\right)$,
respectively, and a number of scalars. The structure of the iteration
step of the method is determined by the envisaged convergence order
$m$ of the method. We consider two strategies for scaling the terms
in the iteration step: norm scaling, where in each step the windows
are normalized, and initial scaling where we only scale in the very
beginning. Norm scaling leads to fast, but conditionally convergent
methods, while initial scaling leads to unconditionally convergent
methods, but with possibly suboptimal convergence constants. The iterations,
initially formulated for time-continuous Gabor systems, are considered
and tested in a discrete setting in which one passes to the appropriately
sampled-and-periodized windows and frame operators. Furthermore, they
are compared with respect to accuracy and efficiency with other methods
to approximate canonical windows associated with Gabor frames.
\end{abstract}
\end{frontmatter}

\section{Introduction}

We consider in this paper iterative schemes for the approximation
of the canonical tight and canonical dual windows associated with
a Gabor frame. We refer to \cite[Ch. 8-10]{Chr03} and \cite[Ch. 5-9, 11-13]{Gro01}
for recent and comprehensive treatments of the theory of Gabor systems
and frames; to fix notations and conventions we briefly give here
the main features. We denote for $g\in L^{2}\left(\mathbb{R}\right)$
and $a>0$, $b>0$ by $\left(g,a,b\right)$ the collection of time-frequency
shifted windows\begin{equation}
g_{na,mb},\quad m,n\in\mathbb{Z},\end{equation}
 where for $x,y\in\mathbb{R}$ we denote\begin{equation}
g_{x,y}=e^{2\pi iyt}g(t-x),\quad t\in\mathbb{R}.\end{equation}
We refer to $g$ as the \emph{window} and to $a$ and $b$ as the
\emph{time-shift} and \emph{frequency-shift} parameters, respectively,
of the Gabor system $\left(g,a,b\right)$. When there are $A>0$,
$B<\infty$, called the \emph{lower} and \emph{upper frame bound},
respectively, such that for all $f\ \in L^{2}\left(\mathbb{R}\right)$
there holds\begin{equation}
A\left\Vert f\right\Vert ^{2}\leq\sum_{m,n=-\infty}^{\infty}\left|\left(f,g_{na,mb}\right)\right|^{2}\leq B\left\Vert f\right\Vert ^{2},\label{eq:frame-inequality}\end{equation}
 we call $\left(g,a,b\right)$ a \emph{Gabor frame}. When in (\ref{eq:frame-inequality})
the second inequality holds for all $f\ \in L^{2}\left(\mathbb{R}\right)$,
we have that \begin{equation}
f\in L^{2}\left(\mathbb{R}\right)\mapsto\, Sf:=\sum_{m,n=-\infty}^{\infty}\left(f,g_{na,mb}\right)g_{na,mb}\end{equation}
 is well-defined as a bounded, positive, semi-definite linear operator
of $L^{2}\left(\mathbb{R}\right)$. We call $S$ the \emph{frame operator}
of $\left(g,a,b\right)$. The frame operator commutes with all relevant
shift operators, i.e., we have for all $f\ \in L^{2}\left(\mathbb{R}\right)$\begin{equation}
Sf_{na,mb}=\left(Sf\right)_{na,mb},\quad m,n\in\mathbb{Z}.\end{equation}

We shall assume in the remainder of this paper that $\left(g,a,b\right)$
is a Gabor frame. Thus the frame operator $S$ is positive definite
and therefore boundedly invertible. There are two windows canonically
associated to the Gabor frame $\left(g,a,b\right)$. These are the
canonical tight window $g^{t}$ and the canonical dual window $g^{d}$,
defined by\begin{eqnarray}
g^{t}=S^{-1/2}g & , & g^{d}=S^{-1}g,\label{eq:gt-gd-def}\end{eqnarray}
respectively. The practical relevance of these windows is that they
give rise to Gabor series representations of arbitrary $f\in L^{2}\left(\mathbb{R}\right)$
according to \begin{equation}
f=\sum_{m,n=-\infty}^{\infty}\left(f,g_{na,mb}^{t}\right)g_{na,mb}^{t}=\sum_{m,n=-\infty}^{\infty}\left(f,g_{na,mb}^{d}\right)g_{na,mb},\end{equation}
where both series are $L^{2}\left(\mathbb{R}\right)$-convergent.
Furthermore, the Gabor systems $\left(g^{t},a,b\right)$ and $\left(g^{d},a,b\right)$
are Gabor frames themselves with frame operators equal to the identity
$I$ and $S^{-1}$, respectively.

The computation of $g^{t}$ and $g^{d}$ according to (\ref{eq:gt-gd-def})
requires taking the inverse square root and the inverse of the frame
operator $S$, respectively. In the often occurring practical case
that $ab$ is a rational number, the frame operator is highly structured
which allows relatively efficient methods for computing $S^{-1}$,
see \cite{qiufei95}. The computation of $S^{-\frac{1}{2}}$ is much
more awkward, even in the case that $ab$ is rational, and often requires
advanced techniques from numerical linear algebra, see for instance
\cite{High99}.

In \cite{jans02} the calculus of Gabor frame operators was combined
with a use of the spectral mapping theorem and Kantorovich's inequality
to analyze an iteration scheme for the approximation of $g^{t}$ that
was proposed around 1995 by Feichtinger and Strohmer (independently
of one another). In this iteration scheme one sets $\gamma_{0}=g$
and for $k=0,1,\ldots$

\begin{minipage}[c]{0.05\columnwidth}%
I.\end{minipage}%
\begin{minipage}[c]{0.95\columnwidth}%
\begin{eqnarray}
\gamma_{k+1} & = & \frac{1}{2}\alpha_{k}\gamma_{k}+\frac{1}{2}\beta_{k}S_{k}^{-1}\gamma_{k};\quad\alpha_{k}=\frac{1}{\left\Vert \gamma_{k}\right\Vert },\,\beta_{k}=\frac{1}{\left\Vert S_{k}^{-1}\gamma_{k}\right\Vert },\label{eq:alg-1}\end{eqnarray}
\end{minipage}%

\noindent where $S_{k}$ is the frame operator corresponding to $\left(\gamma_{k},a,b\right)$.
It was shown in \cite{jans02} that $\left(\gamma_{k},a,b\right)$
is indeed a Gabor frame, and that $\frac{\gamma_{k}}{\left\Vert \gamma_{k}\right\Vert }$
converges (at least) quadratically to $\frac{g^{t}}{\left\Vert g^{t}\right\Vert }$.
From the numerical results in \cite{jans02} for the case that $g$
is the standard Gaussian window $2^{1/4}\exp\left(-\pi t^{2}\right)$
and $a=b=1/\sqrt{2}$ it appears that the resulting method compares
favorably with other iterative techniques for computing inverse square
roots, \cite{kov70,lak95,sher91}.

The investigations in \cite{jans02} were followed by the introduction
in \cite{jans02b,jans03} of two families of iterative algorithms
for the approximation of $g^{t}$ and $g^{d}$, in which the iteration
step involves the initial window $g$ and the frame operator $S$
as well as the current window $\gamma_{k}$ and frame operator $S_{k}$,
but frame operator inversion as in (\ref{eq:alg-1}) do not occur.
The following instances of these two families were analyzed in \cite{jans02b,jans03}.
Again we set $\gamma_{0}=g$, and for $k=0,1,\ldots$

\begin{minipage}[t]{0.05\columnwidth}%
II.\end{minipage}%
\begin{minipage}[c]{0.95\columnwidth}%
\begin{eqnarray}
\gamma_{k+1} & = & \frac{3}{2}\alpha_{k}\gamma_{k}-\frac{1}{2}\beta_{k}S_{k}\gamma_{k};\quad\alpha_{k}=\frac{1}{\left\Vert \gamma_{k}\right\Vert },\,\beta_{k}=\frac{1}{\left\Vert S_{k}\gamma_{k}\right\Vert },\label{eq:alg-2}\end{eqnarray}
\end{minipage}%

\begin{minipage}[t]{0.05\columnwidth}%
III.\end{minipage}%
\begin{minipage}[c]{0.95\columnwidth}%
\begin{eqnarray}
\gamma_{k+1} & = & \frac{15}{8}\varepsilon_{k0}\gamma_{k}-\frac{5}{4}\varepsilon_{k1}S_{k}\gamma_{k}+\frac{3}{8}\varepsilon_{k2}S_{k}^{2}\gamma_{k};\nonumber \\
\varepsilon_{k0} & = & \frac{1}{\left\Vert \gamma_{k}\right\Vert },\,\varepsilon_{k1}=\frac{1}{\left\Vert S_{k}\gamma_{k}\right\Vert },\,\varepsilon_{k2}=\frac{1}{\left\Vert S_{k}^{2}\gamma_{k}\right\Vert },\label{eq:alg-3}\end{eqnarray}
\end{minipage}%

\noindent for the approximation of $g^{t}$, and

\begin{minipage}[t]{0.05\columnwidth}%
IV.\end{minipage}%
\begin{minipage}[c]{0.95\columnwidth}%
\begin{eqnarray}
\gamma_{k+1} & = & 2\alpha_{k}\gamma_{k}-\beta_{k}S_{k}g;\quad\alpha_{k}=\frac{1}{\left\Vert \gamma_{k}\right\Vert },\,\beta_{k}=\frac{1}{\left\Vert S_{k}g\right\Vert },\label{eq:alg-4}\end{eqnarray}
\end{minipage}%

\begin{minipage}[t]{0.05\columnwidth}%
V.\end{minipage}%
\begin{minipage}[c]{0.95\columnwidth}%
\begin{eqnarray}
\gamma_{k+1} & = & 3\delta_{k0}\gamma_{k}-3\delta_{k1}S_{k}g+\delta_{k2}SS_{k}\gamma_{k};\nonumber \\
\delta_{k0} & = & \frac{1}{\left\Vert \gamma_{k}\right\Vert },\,\delta_{k1}=\frac{1}{\left\Vert S_{k}g\right\Vert },\,\delta_{k2}=\frac{1}{\left\Vert SS_{k}\gamma_{k}\right\Vert },\label{eq:alg-5}\end{eqnarray}
\end{minipage}%

\noindent for the approximation of $g^{d}$.

The algorithms II-V are, contrary to algorithm I, conditionally convergent
in the sense that the frame bound ratio $\frac{A}{B}$ of the initial
Gabor frame $\left(g,a,b\right)$ should exceed a certain lower bound.
Accordingly, in algorithm II and III we have that $\frac{\gamma_{k}}{\left\Vert \gamma_{k}\right\Vert }$
converges to $\frac{g^{t}}{\left\Vert g^{t}\right\Vert }$ quadratically
and cubically when $\frac{A}{B}>\frac{1}{2}$ and $\frac{A}{B}>\frac{3}{7}$,
respectively. In algorithm IV and V we have that $\frac{\gamma_{k}}{\left\Vert \gamma_{k}\right\Vert }$
converges to $\frac{g^{d}}{\left\Vert g^{d}\right\Vert }$ quadratically
and cubically when $\frac{A}{B}>\frac{1}{2}\left(\sqrt{5}-1\right)$
and $\frac{A}{B}>0.513829766\ldots$, respectively. A remarkable phenomenon
that emerged from the preliminary experiments done with the algorithms
around 2002, was the fact that the lower bounds for the algorithms
II, III seem far too pessimistic while those for the algorithms IV
and V appear to be realistic.

In the algorithms just presented, all computed windows are normalized.
We shall refer to this as \emph{norm scaling}. Another possibility
that we will investigate, is to replace all scalars $\alpha$'s, $\beta$'s,
$\varepsilon$'s and $\delta$'s that occur in (\ref{eq:alg-2}-\ref{eq:alg-5})
by $1$, and then initially scale the windows by replacing\begin{eqnarray}
g\text{ by }g/\hat{B}^{1/2} & , & S\text{ by }S/\hat{B}.\label{eq:how-to-scale}\end{eqnarray}
 We shall refer to this scaling strategy as \emph{initial scaling}.
If $\hat{B}$ is (an estimate for) the best upper frame bound $\max\sigma\left(S\right)$
then the algorithms will be unconditionally convergent, with guaranteed
desired convergence order, but with convergence constants that may
not be as good as the ones that can be obtained by using the norm
scaling as described by (\ref{eq:alg-2}-\ref{eq:alg-5}).

Matrix versions of algorithms II-III without any scaling have been
treated in \cite{bjobow71}. In \cite{hisc90,kenlaub92} the matrix
version of algorithm II is considered using norm scaling and a scaling
method that approximates the optimal scaling. The matrix version of
algorithm IV is known as a Schulz iteration, see \cite{sch33}. The
fact that $S$, and therefore $\varphi\left(S\right)$ with $\varphi$
continuous and positive on the spectrum of $S$, commutes with all
relevant shift operators, allows us to formulate the iteration steps
on the level of the windows themselves.

In this paper we investigate the algorithms II-V, using both norm
and initial scaling, with more emphasis on computational aspects than
in \cite{jans02,jans02b,jans03}. Here it is necessary to consider
sampled-and-periodized Gabor systems in the style of \cite{Jans97}.
This allows for a formulation and analysis of the algorithms I-V in
an entirely similar way as was done in \cite{jans02,jans02b,jans03}.
Thanks to the fact that the involved (canonical) windows and frame
operators behave so conveniently under the operations of sampling
and periodization, the observations done on the sampled-and-periodized
systems are directly relevant to the time-continuous systems. We must
restrict here to rational values of $ab$, and this gives the frame
operator additional structure which can be exploited in the computations
as dictated by the recursion steps, also see \cite{balfei05,qiufei95,sonder06}
for this matter.

The notions \char`\"{}smart\char`\"{} (or, rather, \char`\"{}smart
but risky\char`\"{}) and \char`\"{}safe\char`\"{} (or, rather, \char`\"{}safe
but conservative\char`\"{}) were introduced in a casual way in \cite{jans03}
to distinguish between cases where the stationary point(s) of the
function transforming (frame) operators according to (\ref{eq:gen-tight-step}),
(\ref{eq:gen-dual-step}) has a good chance to be well-placed in the
middle of and on the \char`\"{}safe\char`\"{} side of the relevant
spectral set, respectively. In the present paper, we choose to refer
to the strategies leading to smart and safe modes as \char`\"{}norm
scaling\char`\"{} and \char`\"{}initial scaling\char`\"{}, respectively,
and discard the terms \char`\"{}smart\char`\"{} and \char`\"{}safe\char`\"{}
altogether.

\section{Paper outline and results}

In Section \ref{sec:Basic-theory} and \ref{sec:Norm-scaling} we
present the basic results of \cite{jans02b,jans03} on transforming
$g$ into $\gamma=\varphi\left(S\right)g$, where $\varphi$ is a
function positive and continuous on $\sigma\left(S\right)$, so as
to obtain a window $\gamma$ whose frame operator $S_{\gamma}$ (for
approximating $g^{t}$) or the operator $\left(SS_{\gamma}\right)^{1/2}=:Z_{\gamma}$
(for approximating $g^{d}$) is closer to (a multiple of) the identity
$I$ than $S$ itself. Here we recall that $\left(g^{t},a,b\right)$
has frame operator $I$ and that $\left(g^{d},a,b\right)$ has frame
operator $S^{-1}$. Thus we relate the operators $S$ and $S_{\gamma}$
and their frame bounds, and we present a bound for the distance between
(the normed) $\gamma$ and $g^{t}$ in terms of the frame bounds of
$S_{\gamma}$. Similarly, we relate the frame bounds of $\left(g,a,b\right)$
and the minimum and maximum of $\sigma\left(Z_{\gamma}\right)$, and
we present a bound for the difference between (the normed) $\gamma$
and $g^{d}$ in terms of the latter minimum and maximum. Next, in
Section \ref{sec:Norm-scaling}, the choice of $\varphi$ is specified
so as to accommodate the recursions of type II, III and of type IV,
V which gives us a means to monitor the frame bound ratio $\frac{A_{k}}{B_{k}}$
(for $g^{t}$) and of the ratio between minimum and maximum of $\sigma\left(Z_{k}\right)$
(for $g^{d}$) during the iteration process.

In Section \ref{sec:Smart-and-safe} we present the algorithms using
only initial scaling.

In Section \ref{sec:Considerations-in-the} we elaborate on the observation
that all algorithms take place in the closed linear span $\mathcal{L}_{g}$
of the adjoint orbit $\left\{ g_{j/b,l/a}\,\big|j,l\in\mathbb{Z}\right\} $.
Here the dual lattice representation of frame operators is relevant
as well as an operator norm to measure the distance of $S_{k}$ (for
$g^{t}$) and of $\left(SS_{k}\right)^{1/2}$ (for $g^{d}$) from
(a multiple of) the identity. The consideration of the algorithms
in the space $\mathcal{L}_{g}$ reveals a fundamental difference between
the algorithms for computing $g^{t}$ and $g^{d}$ that manifests
itself in the totally different after-convergence behaviour of the
two families of algorithms.

In Section \ref{sec:Zak-domain-considerations} we give some considerations
in the Zak transform domain, so as to produce examples of Gabor frames
for which a specific algorithm diverges.

In Section \ref{sec:Sampling-and-periodization} we discuss the discretization
and finitization aspects (through sampling and periodization) that
have to be taken into account since the algorithms are to be tested
numerically.

In Section \ref{sec:Implementational-aspects} we show how the algorithms
can be expressed for discrete, finite Gabor systems, and show that
the algorithms are scalar iterations of the singular values of certain
matrices. We present an efficient implementation of the iterative
algorithms, and we list the window functions we have used to test
the algorithms.

In Section \ref{sec:Experiments} we present our experimental results,
compare them with what the theory predicts and with other methods
to compute tight and dual windows. We provide examples that show the
quadratic and cubic convergence of the algorithms, and the exponential
divergence of the dual iterations after the initial convergence. We
give an example that breaks the norm scaling schemes for both the
tight and dual iterations, and show how various error norms of the
iteration step behave. Comparisons with other methods are made: We
show that the tight iterations are competitive with respect to computing
time and superior with respect to precision. Finally, we show that
the number of iterations needed for full convergence of the algorithms
are dependent on the frame bound ratio, but independent of the structural
properties of the discretization. For initial scaling, we show that
it is easy to choose a scaling parameter that gives almost optimal
convergence.

\section{\label{sec:Basic-theory}Frame operator calculus and basic inequalities}

The basic theory to analyze the recursions appears somewhat scattered
in \cite{jans02,jans02b,jans03}; for the reader's convenience, we
give in Section \ref{sec:Basic-theory} and \ref{sec:Norm-scaling}
a concise yet comprehensive summary of the basic results and ideas.
We let $\left(g,a,b\right)$ be a Gabor frame with frame operator
$S$ and best frame bounds $A=\min\sigma\left(S\right)>0$, $B=\max\sigma\left(S\right)$,
where $\sigma\left(S\right)$ is the spectrum of $S$. In this section
we present the basic inequalities expressing the approximation errors
in terms of the (frame) bounds on the involved (frame) operators.
These inequalities are a consequence of the calculus of Gabor frame
operators, the spectral mapping theorem and Kantorovich's inequality.

\begin{prop}
Let $\varphi$ be continuous and positive on $\left[A,B\right]$,
and set $\gamma:=\varphi\left(S\right)g$. The following holds.
\begin{lyxlist}{(ii)}
\item [(i)]$\left(\gamma,a,b\right)$ is a Gabor frame with frame operator
$S_{\gamma}:=S\varphi^{2}\left(S\right)$ and best frame bounds\begin{eqnarray}
A_{\gamma}:=\min_{s\in\sigma\left(S\right)}s\varphi^{2}\left(s\right) & , & \, B_{\gamma}:=\max_{s\in\sigma\left(S\right)}s\varphi^{2}\left(s\right).\end{eqnarray}
Furthermore, \begin{equation}
g^{t}=S^{-1/2}g=S_{\gamma}^{-1/2}\gamma=\gamma^{t},\end{equation}
and\begin{equation}
\left\Vert \frac{\gamma}{\left\Vert \gamma\right\Vert }-\frac{g^{t}}{\left\Vert g^{t}\right\Vert }\right\Vert \leq\left(1-Q_{\gamma}^{1/4}\right)\sqrt{\frac{2}{1+Q_{\gamma}}};\quad Q_{\gamma}=\frac{A_{\gamma}}{B_{\gamma}}.\label{eq:kant-t}\end{equation}

\item [(ii)]Let $Z_{\gamma}:=\left(SS_{\gamma}\right)^{1/2}=S\varphi\left(S\right)$,
and \begin{eqnarray}
E_{\gamma} & := & \min\sigma\left(Z_{\gamma}\right)=\min_{s\in\sigma\left(S\right)}s\varphi\left(s\right),\\
F_{\gamma} & := & \max\sigma\left(Z_{\gamma}\right)=\max_{s\in\sigma\left(S\right)}s\varphi\left(s\right).\end{eqnarray}
Then\begin{eqnarray}
g^{d}=Z_{\gamma}^{-1}\gamma & , & Z_{\gamma}\gamma=S_{\gamma}g,\label{eq:funinv-d}\end{eqnarray}
and\begin{equation}
\left\Vert \frac{\gamma}{\left\Vert \gamma\right\Vert }-\frac{g^{d}}{\left\Vert g^{d}\right\Vert }\right\Vert \leq\left(1-R_{\gamma}^{1/2}\right)\sqrt{\frac{2}{1+R_{\gamma}}};\quad R_{\gamma}=\frac{E_{\gamma}}{F_{\gamma}}.\label{eq:kant-d}\end{equation}

\end{lyxlist}
\end{prop}
The basic result (i) gives us a clue how to produce a good approximation
$\frac{\gamma}{\left\Vert \gamma\right\Vert }$ of $\frac{g^{t}}{\left\Vert g^{t}\right\Vert }$:
Take $\varphi$ such that $s\varphi^{2}\left(s\right)$ is flat on
$\sigma\left(S\right)\subset\left[A,B\right]$ so that the number
$Q_{\gamma}$ in (\ref{eq:kant-t}) is close to $1$. Similarly, by
the basic result (ii), the number $R_{\gamma}$ in (\ref{eq:kant-d})
is close to $1$ when $\varphi$ is such that $s\varphi\left(s\right)$
is flat on $\sigma\left(S\right)\subset\left[A,B\right]$, and then
we obtain a good approximation of $\frac{g^{d}}{\left\Vert g^{d}\right\Vert }$.
In the next two subsections, we use this basic result repeatedly with
polynomials $\varphi$ of fixed degree $m$ so as to obtain iterative
approximations of $g^{t}$ and $g^{d}$.

\section{Norm scaling\label{sec:Norm-scaling}}

\subsection{\label{sub:Iterations-for-t}Iterations for approximating $g^{t}$}

We consider iteration schemes\begin{eqnarray}
\gamma_{0} & = & g;\quad\gamma_{k+1}=\varphi_{k}\left(S_{k}\right)\gamma_{k},\quad k=0,1,\ldots,\end{eqnarray}
for the approximation of $g^{t}$, where $S_{k}$ is the frame operator
of $\left(\gamma_{k},a,b\right)$. We use here the basic result (i)
repeatedly with\begin{eqnarray}
g & = & \gamma_{k},\quad S=S_{k}\textrm{ and }\gamma=\gamma_{k+1},\quad S_{\gamma}=S_{k+1}=S_{k}\varphi_{k}^{2}\left(S_{k}\right).\label{eq:gen-tight-step}\end{eqnarray}
 For $k=0,1,\ldots$ we have that $\gamma_{k}^{t}=g^{t}$ and that\begin{equation}
\left\Vert \frac{\gamma_{k}}{\left\Vert \gamma_{k}\right\Vert }-\frac{g^{t}}{\left\Vert g^{t}\right\Vert }\right\Vert \leq\left(1-Q_{k}^{1/4}\right)\sqrt{\frac{2}{1+Q_{k}}};\quad Q_{k}=\frac{A_{k}}{B_{k}},\end{equation}
where $A_{k}$ and $B_{k}$ are the best frame bounds of $S_{k}$.
The numbers $A_{k}$, $B_{k}$ can be computed and estimated recursively
according to $A_{0}=A$, $B_{0}=B$ and \begin{eqnarray}
A_{k+1} & = & \min_{s\in\sigma\left(S_{k}\right)}s\varphi_{k}^{2}\left(s\right)\geq\min_{s\in\left[A_{k},B_{k}\right]}s\varphi_{k}^{2}\left(s\right),\\
B_{k+1} & = & \max_{s\in\sigma\left(S_{k}\right)}s\varphi_{k}^{2}\left(s\right)\leq\max_{s\in\left[A_{k},B_{k}\right]}s\varphi_{k}^{2}\left(s\right),\end{eqnarray}
for $k=0,1,\ldots$.

We should choose $\varphi_{k}$ such that $s\varphi_{k}^{2}\left(s\right)$
is flat on $\sigma\left(S_{k}\right)$. To that end there is proposed
in \cite[Subsec. 5.1]{jans03} for $m=2,3,\ldots$ the choice\begin{equation}
\varphi_{k}\left(s\right)=\sum_{j=0}^{m-1}a_{mj}\alpha_{kj}s^{j},\quad\alpha_{kj}=\left\Vert S_{k}^{j}\gamma_{k}\right\Vert ^{-1},\label{eq:phi-t-as-pol}\end{equation}
where the $a_{mj}$ are defined by\begin{equation}
\sum_{l=0}^{m-1}\left(-1\right)^{l}\left(\begin{array}{c}
-1/2\\
l\end{array}\right)\left(1-x\right)^{l}=\sum_{j=0}^{m-1}a_{mj}x^{j},\quad x>0.\label{eq:a_kj}\end{equation}
The motivation for this choice is as follows. The left-hand side of
(\ref{eq:a_kj}) is the $\left(m-1\right)^{\textrm{th}}$ order Taylor
approximation of $x^{-1/2}$ around $x=1$, while (when $Q_{k}$ is
sufficiently close to 1)\begin{eqnarray}
\alpha_{kj}S_{k}^{j} & \approx & \left(\frac{S_{k}}{\left\Vert S_{k}\right\Vert }\right)^{j}\frac{1}{\left\Vert \gamma_{k}\right\Vert },\quad j=0,\ldots,m-1.\end{eqnarray}
Hence $s\varphi_{k}^{2}\left(s\right)=\left(s^{1/2}\varphi_{k}\left(s\right)\right)^{2}$
should be expected to be flat on $\sigma\left(S_{k}\right)$, with
$1-Q_{k+1}$ potentially of order $\left(1-Q_{k}\right)^{m}$.

When $m=2,3$ we get the iterations II, III in (\ref{eq:alg-2}),
(\ref{eq:alg-3}). It is shown in \cite[Sec. 4]{jans02b} and \cite[Sec. 6]{jans03}
that for $m=2$ the quantity $Q_{k}=\frac{A_{k}}{B_{k}}$ increases
to $1$ and that $\frac{\gamma_{k}}{\left\Vert \gamma_{k}\right\Vert }\rightarrow\frac{g^{t}}{\left\Vert g^{t}\right\Vert }$
quadratically when $k\rightarrow\infty$, provided that $\frac{A}{B}>\frac{1}{2}$.
For $m=3$ it is shown in \cite[Sec. 8]{jans02b} that $Q_{k}$ increases
to $1$ and that $\frac{\gamma_{k}}{\left\Vert \gamma_{k}\right\Vert }\rightarrow\frac{g^{t}}{\left\Vert g^{t}\right\Vert }$
cubically when $k\rightarrow\infty$, provided that $\frac{A}{B}>\frac{3}{7}$.
In \cite{jans02b,jans03} it was observed for $m=2,3$ that the choice
of $\alpha_{kj}$ causes $s\varphi_{k}^{2}\left(s\right)$ to have
one or more stationary points in $\left[A_{k},B_{k}\right]$ so that
the odds for flatness of $s\varphi_{k}^{2}\left(s\right)$ on $\left[A_{k},B_{k}\right]$
are favourable.

\subsection{\label{sub:Iterations-for-d}Iterations for approximating $g^{d}$}

We consider iteration schemes\begin{eqnarray}
\gamma_{0} & = & g;\quad\gamma_{k+1}=\varphi_{k}\left(Z_{k}\right)\gamma_{k},\quad k=0,1,\ldots,\label{eq:gen-dual-step}\end{eqnarray}
for the approximation of $g^{d}$, where $Z_{k}=\left(SS_{k}\right)^{1/2}$
with $S_{k}$ the frame operator of $\left(\gamma_{k},a,b\right)$.
It is seen from (\ref{eq:gen-dual-step}) by induction that $\gamma_{k}=\psi_{k}\left(S\right)g$
for some function $\psi_{k}$. Hence by the basic result (i) with
$\varphi=\psi_{k}$, we have that\begin{eqnarray}
S_{k}=S\psi_{k}^{2}\left(S\right) & , & Z_{k}=S\psi_{k}\left(S\right),\end{eqnarray}
 and, by the basic result (ii), that \begin{equation}
\left\Vert \frac{\gamma_{k}}{\left\Vert \gamma_{k}\right\Vert }-\frac{g^{d}}{\left\Vert g^{d}\right\Vert }\right\Vert \leq\left(1-R_{k}^{1/2}\right)\sqrt{\frac{2}{1+R_{k}}};\quad R_{k}=\frac{E_{k}}{E_{k}},\end{equation}
 where $E_{k}=\min\sigma\left(Z_{k}\right)$, $F_{k}=\max\sigma\left(Z_{k}\right)$.
A further use of the calculus of frame operators as given by the basic
result (i), yields\begin{equation}
Z_{k+1}=Z_{k}\varphi_{k}\left(Z_{k}\right).\end{equation}
Consequently, the numbers $E_{k}$, $F_{k}$ can be computed and estimated
recursively according to $E_{0}=A$, $F_{0}=B$ and\begin{eqnarray}
E_{k+1} & = & \min_{z\in\sigma\left(Z_{k}\right)}z\varphi_{k}\left(z\right)\geq\min_{z\in\left[E_{k},F_{k}\right]}z\varphi_{k}\left(z\right),\\
F_{k+1} & = & \max_{z\in\sigma\left(Z_{k}\right)}z\varphi_{k}\left(z\right)\leq\max_{z\in\left[E_{k},F_{k}\right]}z\varphi_{k}\left(z\right),\end{eqnarray}
for $k=0,1,\ldots$.

We should choose $\varphi_{k}$ such that $z\varphi_{k}\left(z\right)$
is flat on $\sigma\left(Z_{k}\right)$. To that end there is proposed
in \cite[Subsec. 5.2]{jans03} for $m=2,3,\ldots$ the choice\begin{equation}
\varphi_{k}\left(z\right)=\sum_{j=0}^{m-1}b_{mj}\beta_{kj}z^{j},\quad\beta_{kj}=\left\Vert Z_{k}^{j}\gamma_{k}\right\Vert ^{-1},\end{equation}
where the $b_{mj}$ are defined by\begin{equation}
\sum_{l=0}^{m-1}\left(1-x\right)^{l}=\sum_{j=0}^{m-1}b_{mj}x^{j},\quad x>0.\label{eq:b_kj}\end{equation}
The motivation for the proposal is similar to the one for the choice
of $\varphi_{k}$ in (\ref{eq:phi-t-as-pol}) in Subsec. \ref{sub:Iterations-for-t};
we now note that the left-hand side of (\ref{eq:b_kj}) is the $\left(m-1\right)^{\textrm{th}}$
order Taylor approximation of $x^{-1}$ around $x=1$. The implementation
of the resulting recurrence step\begin{eqnarray}
\gamma_{k+1} & = & \sum_{j=0}^{m-1}b_{mj}\frac{Z_{k}^{j}\gamma_{k}}{\left\Vert Z_{k}^{j}\gamma_{k}\right\Vert },\quad Z_{k}=\left(SS_{k}\right)^{1/2},\end{eqnarray}
is made feasible by the observation that, thanks to the second item
in (\ref{eq:funinv-d}), $Z_{k}\gamma_{k}=S_{k}g$ so that\begin{eqnarray}
Z_{k}^{2r}\gamma_{k} & = & \left(SS_{k}\right)^{r}\gamma_{k},\quad Z_{k}^{2r+1}\gamma_{k}=\left(SS_{k}\right)^{r}S_{k}g,\quad r=0,1,\ldots.\end{eqnarray}

When $m=2,3$ we get the recurrences IV, V in (\ref{eq:alg-4}) and
(\ref{eq:alg-5}). It is shown in \cite[Sec. 5]{jans02b} and \cite[Sec. 7]{jans03},
that for $m=2$ the quantity $R_{k}=\frac{E_{k}}{F_{k}}$ increases
to $1$ and that $\frac{\gamma_{k}}{\left\Vert \gamma_{k}\right\Vert }\rightarrow\frac{g^{d}}{\left\Vert g^{d}\right\Vert }$
quadratically when $k\rightarrow\infty$, provided that $\frac{A}{B}>\frac{1}{2}\left(\sqrt{5}-1\right)$.
For $m=3$ it is shown in \cite[Sec. 9]{jans02b} that $R_{k}$ increases
to $1$ and that $\frac{\gamma_{k}}{\left\Vert \gamma_{k}\right\Vert }\rightarrow\frac{g^{d}}{\left\Vert g^{d}\right\Vert }$
cubically when $k\rightarrow\infty$, provided that $\frac{A}{B}>0.513829766\ldots$.
In \cite{jans02b,jans03} it was observed for $m=2,3$ that the choice
of $\beta_{kj}$ causes $z\varphi_{k}\left(z\right)$ to have one
or more stationary points in $\left[E_{k},F_{k}\right]$.

\section{\label{sec:Smart-and-safe}Initial scaling}

The algorithms II-V are guaranteed to converge when the lower bound
ratio $\frac{A}{B}$ of $\left(g,a,b\right)$ exceeds a certain value.
The proofs, as given in \cite{jans02b} and \cite{jans03}, require
a careful analysis of the extreme values of the functions $\varphi_{k}$
on the spectra of the relevant operators and can become quite complicated,
especially in the cases of algorithms III, V. However, the algorithms
are efficient in the sense that the envisaged convergence order $m$
is realized with favourable convergence constants. In practice, as
the experiments in Section \ref{sec:Experiments} show, the algorithms
II, III turn out to converge in almost all cases, even when the frame
bound ratio is close to $0$. However, divergence of the algorithms
IV, V occurs much more frequently. In Section \ref{sec:Zak-domain-considerations}
we present examples, using the Zak transform, of frames $\left(g,a,b\right)$
such that algorithm II and IV diverges.

It would be desirable to have versions of the algorithms that are
guaranteed to converge, no matter how small the frame bound ratio
of the initial frame is (as long as it is positive). In the following,
we present the initial scaling versions of the algorithms that converge
at the envisaged convergence order $m$, possibly with suboptimal
convergence constants. Since we can freely switch scaling strategy,
a possible strategy is to initially scale such that convergence is
guaranteed, and to continue until one is confident that the relevant
condition number exceeds the specific lower bound so that the norm
scaling mode can be applied from that point onwards.

The introduction in \cite{jans03} of the notion of {}``safe modes''
was prompted by an observation by M. Hampejs who prescaled the window
$g$ (and the frame operator) and deleted all normalization operations
in the recursion step of algorithms II, IV. In the present paper,
the prescaling is done in such a way that the scaled $S$ has its
spectrum exclusively in the attraction region of the function $\varphi$
describing the simplified recursion. More specifically, we consider
the iteration steps as given in Subsections \ref{sub:Iterations-for-t},
\ref{sub:Iterations-for-d}, with all $\alpha$'s and $\beta$'s equal
to $1$. The $\varphi$'s thus obtained are independent of $k$ and
are given by\begin{eqnarray}
\varphi_{m}^{t}\left(s\right) & := & \sum_{l=0}^{m-1}\left(-1\right)^{l}\left(\begin{array}{c}
-1/2\\
l\end{array}\right)\left(1-s\right)^{l}=\sum_{j=0}^{m-1}a_{mj}s^{j},\quad s>0,\label{eq:phi_m_t}\end{eqnarray}
and\begin{eqnarray}
\varphi_{m}^{d}\left(z\right) & := & \sum_{l=0}^{m-1}\left(1-z\right)^{l}=\sum_{j=0}^{m-1}b_{mj}z^{j},\quad z>0,\label{eq:phi_m_d}\end{eqnarray}
respectively. The relevant spectra transform by the spectral mapping
theorem according to\begin{eqnarray}
\sigma\left(S_{k}\right) & \rightarrow & \left\{ s\left(\varphi_{m}^{t}\left(s\right)\right)^{2}\,\bigg|\, s\in\sigma\left(S_{k}\right)\right\} =\sigma\left(S_{k+1}\right),\end{eqnarray}
and\begin{eqnarray}
\sigma\left(Z_{k}\right) & \rightarrow & \left\{ z\varphi_{m}^{d}\left(z\right)\,\bigg|\, z\in\sigma\left(Z_{k}\right)\right\} =\sigma\left(Z_{k+1}\right),\end{eqnarray}
respectively.

\begin{figure}
\subfigure[$s>0\mapsto s\left(\varphi_{m}^{t}\left(s\right)\right)^{2}$.]{\includegraphics[%
  width=0.5\textwidth,
  keepaspectratio]{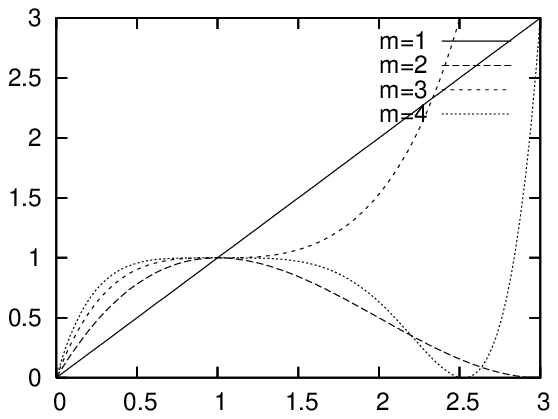}}\subfigure[$z>0\mapsto z\varphi_{m}^{d}\left(z\right).$]{

\includegraphics[%
  width=0.5\textwidth,
  keepaspectratio]{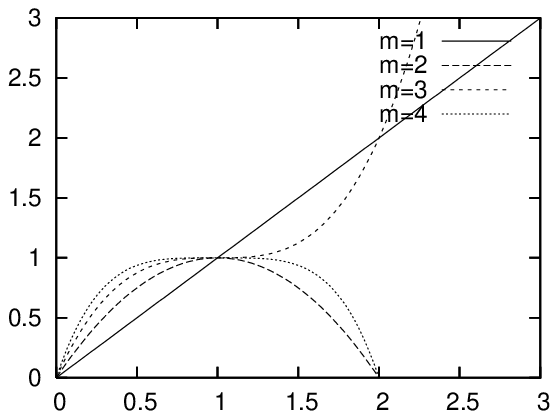}}

\caption{\label{cap:phi}The figure shows the two set of functions governing
the convergence of the tight and dual iterations using initial scaling
for order $m=1,2,3,4$ with $\varphi_{m}^{t}$ and $\varphi_{m}^{d}$
defined by (\ref{eq:phi_m_t}) and (\ref{eq:phi_m_d}). For (a), the
attraction point $1$ has attraction regions $m=2:\,\left(0,3\right)$,
$m=3:\,\left(0,7/3\right)$, $m=4:\,\left(0,2.525847988\right)$.
For (b), the attraction point $1$ has attraction region $(0,2)$
for $m=2,3,4$. }
\end{figure}

The functions $\varphi_{m}^{t}$ and $\varphi_{m}^{d}$ are $\left(m-1\right)^{\textrm{th}}$
order Taylor approximations of $s^{-1/2}$ and $z^{-1}$ around $s=1$
and $z=1$, respectively. Hence $s\left(\varphi_{m}^{t}\left(s\right)\right)^{2}$
and $z\varphi_{m}^{d}\left(z\right)$ approximate $1$ around $s=1$
and $z=1$, respectively. In Fig. \ref{cap:phi} we have shown plots
of the mappings\begin{eqnarray}
s>0\mapsto s\left(\varphi_{m}^{t}\left(s\right)\right)^{2} & , & z>0\mapsto z\varphi_{m}^{d}\left(z\right)\end{eqnarray}
 for $m=1,2,3,4$, respectively. Fig. \ref{cap:phi}(a) also appears
in \cite{bjobow71}. In all cases, the point $s=1=z$ is an attractor
for the region $\left(0,2\right)$. Consequently, when $S_{0}=S$,
$Z_{0}=S$ have spectrum in $\left(0,2\right)$, the spectra $\sigma\left(S_{k}\right)$,
$\sigma\left(Z_{k}\right)$ converge to $1$ as $k\rightarrow\infty$,
and the convergence is of order $m$ in the sense that the ratio of
minimum and maximum of the spectra converge to $1$ at order $m$.
Thus we should replace $g$ by $g/\left(\hat{B}\right)^{1/2}$ and
$S$ by $S/\hat{B}$ where $\hat{B}$ is such that $\sigma\left(S/\hat{B}\right)\subset\left(0,2\right)$
to obtain iterations having $m^{\textrm{th}}$ order convergence to
$g^{t}$ and to $\left(\hat{B}\right)^{1/2}g^{d}$, respectively.

To guarantee convergence an estimate of $\max\sigma\left(S\right)$
is needed. In \cite{balfei05} a number of upper bounds of $\max\sigma\left(S\right)$
are developed for discrete-time, periodic Gabor systems. A convenient
upper bound for our purposes follows from the dual lattice representation
of the frame operator $S$, see Section \ref{sec:Considerations-in-the}
for more details, as\begin{equation}
\max\sigma\left(S\right)\leq\frac{1}{ab}\sum_{j,l}\left|\left(g,g_{j/b,l/a}\right)\right|.\end{equation}

We make some comments for scaling optimally in the first iteration
step. We shall refer to this method as \emph{initial optimal scaling}.
Assume that $\sigma\left(S\right)$ consists of the entire interval
$[A,B]$. Consider the tight iterations as described in this Section,
and assume that we replace $S$ by $S/\hat{B}$. Then\begin{eqnarray}
A_{1} & = & \min\left\{ s\left(\varphi_{m}^{t}\left(s\right)\right)^{2}\bigg|s\in\left[A/\hat{B},B/\hat{B}\right]\right\} ,\\
B_{1} & = & \max\left\{ s\left(\varphi_{m}^{t}\left(s\right)\right)^{2}\bigg|s\in\left[A/\hat{B},B/\hat{B}\right]\right\} .\end{eqnarray}

\begin{table}

\caption{\label{cap:optscal}}

\begin{longtable}{r|c}
\hline 
Method.&
$\hat{B}$ or $\hat{F}$\tabularnewline
\hline
I.&
$\sqrt{AB}$\tabularnewline
II.&
$\frac{1}{3}\left(A+\sqrt{AB}+B\right)$\tabularnewline
III.&
$\frac{3}{10}\left(B+A\right)+\frac{2}{5}\sqrt{\frac{1}{2}\left(B^{2}+A^{2}\right)+\frac{1}{16}\left(B-A\right)^{2}}$\tabularnewline
IV.&
$\frac{1}{2}\left(E+F\right)$\tabularnewline
V.&
$\frac{1}{3}\left(F+E\right)+\frac{1}{3}\sqrt{\frac{1}{2}\left(F^{2}+E^{2}\right)+\frac{1}{2}\left(F-E\right)^{2}}$\tabularnewline
\hline
\end{longtable}

The table shows the optimal scaling constant, $\hat{B}$ or $\hat{F}$,
for doing initial scaling of the five iteration types.\\

\end{table}

Initial optimal scaling occurs for that value of $\hat{B}$ for which
the ratio $A_{1}/B_{1}$ is maximal. The optimal value of the scaling
parameter $\hat{F}$ for the dual iterations is defined in a similar
way. For the five iteration types, the optimal $\hat{B}$ or $\hat{F}$
is shown in Table \ref{cap:optscal}. All these numbers are close
to the center of the interval $\left[A,B\right]$ or $\left[E,F\right]$.
It is not hard to show that all optimally scaled operators $S$ and
$Z$ have their spectra in the attraction regions given in Figure
\ref{cap:phi} for the algorithms II-V.

If we scale optimally in each iteration step, and not just the first,
we get the best possible convergence constants. However, this method
is not practically feasible because of the repeated calculations of
frame bounds, and we shall use it only as a reference method. We refer
to it as \emph{constant optimal scaling}, see Fig. \ref{cap:conv_scal}.

\section{\label{sec:Considerations-in-the}Considerations in the adjoint orbit
space}

We consider the closed linear span $\mathcal{L}_{g}$ of the adjoint
orbit $\left\{ g_{j/b,l/a}\,\big|\, j,l\in\mathbb{Z}\right\} $. According
to the duality principle of Gabor analysis we have that the adjoint
orbit is a Riesz basis for $\mathcal{L}_{g}$ (for this matter we
refer to \cite[Secs. 3.6 and 9.2]{Chr03}, and \cite[Ch. 7]{Gro01}).
Furthermore, when $f\in L^{2}\left(\mathbb{R}\right)$, the orthogonal
projection of $f$ onto $\mathcal{L}_{g}$ is given by 

\begin{eqnarray}
P_{g}f & = & \frac{1}{ab}\sum_{j,l}\left(f,\left(g^{d}\right)_{j/b,l/a}\right)g_{j/b,l/a}=\frac{1}{ab}\sum_{j,l}\left(f,\left(g^{t}\right)_{j/b,l/a}\right)\left(g^{t}\right)_{j/b,l/a}.\end{eqnarray}
As a consequence of $\gamma_{k}^{t}=g^{t}$ in all our algorithms,
we see that $\gamma_{k}\in\mathcal{L}_{g}$. There is also the Wexler-Raz
biorthogonality relation,\begin{eqnarray}
\left(g,\left(g^{d}\right)_{j/b,l/a}\right) & = & \left(g^{t},\left(g^{t}\right)_{j/b,l/a}\right)=ab\delta_{j0}\delta_{l0},\label{eq:wexraz}\end{eqnarray}
where $\delta$ is Kronecker's delta. Finally, there is the following
fundamental identity of Gabor analysis. Assume that $f,\xi,\gamma,h\in L^{2}\left(\mathbb{R}\right)$
and that the three Gabor systems $\left(f,a,b\right)$, $\left(\xi,a,b\right)$,
$\left(\gamma,a,b\right)$ have finite upper frame bounds. Then we
have\begin{eqnarray}
\sum_{m,n}\left(f,\gamma_{na,mb}\right)\left(\xi_{na,mb},h\right) & = & \frac{1}{ab}\sum_{j,l}\left(\xi,\gamma_{j/b,l/a}\right)\left(f_{j/b,l/a},h\right)\label{eq:fund-ident}\end{eqnarray}
with absolute convergence at either side. We can regard (\ref{eq:fund-ident})
as a representation result for the frame-type operator\begin{eqnarray}
S_{\gamma,\xi} & : & f\rightarrow\sum_{n,m}\left(f,\gamma_{na,mb}\right)\xi_{na,mb},\end{eqnarray}
viz. as\begin{eqnarray}
S_{\gamma,\xi} & = & \frac{1}{ab}\sum_{j,l}\left(\xi,\gamma_{j/b,l/a}\right)U_{j,l},\label{eq:jans-rep}\end{eqnarray}
where $U_{j,l}$ is the unitary operator\begin{eqnarray}
U_{j,l} & : & f\rightarrow f_{j/b,l/a}.\end{eqnarray}
This is the dual lattice representation (also known as the Janssen
representation, see \cite[Sec. 7.2]{Chr03} and \cite[Corr. 9.3.7]{Gro01})
of the frame operator. In order for (\ref{eq:jans-rep}) to be well-defined,
we assume that $\gamma,\xi$ satisfies the so-called \emph{Condition
A}':

\begin{lyxlist}{MMM.}
\item [A':]\begin{equation}
\sum_{j,l}\left|\left(\xi,\gamma_{j/b,l/a}\right)\right|<\infty,\end{equation}

\end{lyxlist}
see \cite[Def. 7.2.1]{Gro01}. If $\xi=\gamma$ then this is the Condition
A introduced by Tolimieri and Orr in \cite{tolorr95}. We refer to
Appendix \ref{sec:Condition-A} where an instance, relevant in the
present context, of a pair $\xi,\gamma$ satisfying condition A' is
given.

\subsection{Estimate for upper frame bound}

If $g$ satisfies Condition A, the frame operator $S$ of $\left(g,a,b\right)$
has the representation\begin{eqnarray}
S & = & \frac{1}{ab}\sum_{j,l}\left(g,g_{j/b,l/a}\right)U_{j,l},\end{eqnarray}
with absolute convergence in the operator norm. Therefore, there is
the upper bound\begin{eqnarray}
\hat{B} & = & \frac{1}{ab}\sum_{j,l}\left|\left(g,g_{j/b,l/a}\right)\right|\label{eq:Bbound}\end{eqnarray}
 for the best upper frame bound $\max\sigma\left(S\right)$ of $\left(g,a,b\right)$.

\subsection{Error measure}

We measure convergence of $\gamma_{k}$ to $g^{t}$ and $g^{d}$ by
inspecting $L^{2}$-distances of the normed windows. This quantity
is bounded in terms of the numbers $Q_{k}$ and $R_{k}$ in (\ref{eq:kant-t})
and (\ref{eq:kant-d}) that measure how close the operators $S_{k}$
and $Z_{k}$ are to being a multiple of the identity operator. In
the converse direction, it would be useful to have a measure on the
windows that translates directly to the distance of $S_{k}$ and $Z_{k}$
to (a multiple of) the identity operator. Such a measure can indeed
be found. As to $g^{t}$ we note that when $g$ satisfies Condition
A then $S_{k}$ has the representation\begin{eqnarray}
S_{k} & = & \frac{1}{ab}\sum_{j,l}\left(\gamma_{k},\left(\gamma_{k}\right)_{j/b,l/a}\right)U_{j,l},\end{eqnarray}
whence\begin{eqnarray}
\left\Vert S_{k}-\frac{1}{ab}\left\Vert \gamma_{k}\right\Vert ^{2}I\right\Vert  & \leq & \frac{1}{ab}\sum_{j,l\neq0,0}\left|\left(\gamma_{k},\left(\gamma_{k}\right)_{j/b,l/a}\right)\right|.\label{eq:t-jans-bound}\end{eqnarray}
As to $g^{d}$ we note that $Z_{k}=S\varphi_{k}\left(S\right)$ and
with $\gamma_{k}=\varphi_{k}\left(S\right)g$ there holds by frame
operator calculus\begin{eqnarray}
S\varphi_{k}\left(S\right)f & = & \sum_{m,n}\left(f,\left(\gamma_{k}\right)_{na,mb}\right)g_{na,mb}.\end{eqnarray}
Hence there is the representation\begin{eqnarray}
Z_{k} & = & S\varphi_{k}\left(S\right)=\frac{1}{ab}\sum_{j,l}\left(g,\left(\gamma_{k}\right)_{j/b,l/a}\right)U_{j,l}.\end{eqnarray}
Therefore\begin{equation}
\left\Vert Z_{k}-\frac{1}{ab}\left(g,\gamma_{k}\right)I\right\Vert \leq\frac{1}{ab}\sum_{j,l\neq0,0}\left|\left(g,\left(\gamma_{k}\right)_{j/b,l/a}\right)\right|.\label{eq:d-jans-bound}\end{equation}
 Note that the quantities of the right-hand sides of (\ref{eq:t-jans-bound})
and (\ref{eq:d-jans-bound}) measure to what extent the Wexler-Raz
condition (\ref{eq:wexraz}) is violated. In Appendix \ref{sec:Condition-A}
it is shown for $a\in\mathbb{N}$, $b^{-1}\in\mathbb{N}$ and $g$
satisfying condition A that the $\gamma_{k}$ occurring in (\ref{eq:t-jans-bound})
and the $g$, $\gamma_{k}$ occurring in (\ref{eq:d-jans-bound})
satisfy condition A and A', respectively. We shall refer to the right
hand sides of (\ref{eq:t-jans-bound}) and (\ref{eq:d-jans-bound})
as the \emph{dual lattice norm}.

\subsection{\label{sub:Influence-of-out-of-space}Influence of out-of-space components}

We have seen that all iterands $\gamma_{k}$ of the algorithms are
in $\mathcal{L}_{g}$. We briefly comment on the impact on the algorithms
of $\gamma_{k}$ having non-zero components orthogonal to $\mathcal{L}_{g}$
(one can think here of round-off errors generating these components).
To that end we consider the algorithms II and IV (assuming appropriate
scaling has been carried out), and we assume that they have converged
to the extent that the operators $S_{k}$ and $Z_{k}$ agree within
machine precision with the identity operator.

As for algorithm II we thus have that\begin{eqnarray}
\gamma_{k+1} & = & \frac{3}{2}\gamma_{k}-\frac{1}{2}S_{k}\gamma_{k}=\gamma_{k}\end{eqnarray}
 within machine precision. Hence, possible out-of-space components
in $\gamma_{k}$ are reproduced within machine precision. As a consequence,
we should expect that the error stays at its converged level when
the iteration is continued beyond the point where machine precision
is reached.

Next we consider algorithm IV using initial scaling so that\begin{eqnarray}
\gamma_{k+1} & = & 2\gamma_{k}-S_{k}g.\label{eq:IV-error}\end{eqnarray}
The term $S_{k}g$ has the representation \begin{eqnarray}
S_{k}g & = & \frac{1}{ab}\sum_{j,l}\left(\gamma_{k},\left(\gamma_{k}\right)_{j/b,l/a}\right)g_{j/b,l/a}\in\mathcal{L}_{g}.\end{eqnarray}
Furthermore,\begin{eqnarray}
P_{g}\gamma_{k} & = & \frac{1}{ab}\sum_{j,l}\left(\gamma_{k},\left(g^{d}\right)_{j/b,l/a}\right)g_{j/b,l/a}.\end{eqnarray}
Hence, $S_{k}g=P_{g}\gamma_{k}\in\mathcal{L}g$ to machine precision,
and, to machine precision,\begin{eqnarray}
P_{g}\gamma_{k+1} & = & \frac{1}{ab}\sum_{j,l}\left(\gamma_{k},\left(\gamma_{k}\right)_{j/b,l/a}\right)g_{j/b,l/a}=P_{g}\gamma_{k}.\end{eqnarray}
On the other hand the orthogonal component $\gamma_{k}-P_{g}\gamma_{k}$
is per (\ref{eq:IV-error}) multiplied by $2$, i.e., to machine precision,\begin{eqnarray}
\gamma_{k+1}-P_{g}\gamma_{k+1} & = & 2\left(\gamma_{k}-P_{g}\gamma_{k}\right).\end{eqnarray}
 As a consequence, the algorithm starts to diverge beyond the point
where machine precision is reached.

The observations just made continue to hold for the more general algorithms
in Subsections \ref{sub:Iterations-for-t} and \ref{sub:Iterations-for-d}.
Thus no substantial after-convergence error build-up occurs for the
algorithms of Subsection \ref{sub:Iterations-for-t}. For the algorithms
of Subsection \ref{sub:Iterations-for-d}, with basic recursion step\begin{eqnarray}
\gamma_{k+1} & = & \sum_{j=0}^{m-1}b_{mj}Z_{k}^{j}\gamma_{k}\end{eqnarray}
 the terms with odd $j$ all lie in $\mathcal{L}_{g}$, and those
with even $j$ are given within machine precision by $b_{mj}\gamma_{k}$.
Since $\sum_{j\textrm{ even}}b_{mj}=2^{m-1}$, the out-of-space component
in $\gamma_{k}$ gets multiplied by $2^{m-1}$ in each iteration step.

\section{\label{sec:Zak-domain-considerations}Zak domain considerations}

We consider the case that $ab=$$\frac{p}{q}$ with integer $p,q>0$
such that $\gcd\left(p,q\right)=1$, and we define the Zak transform
$Z$ as (the extension to $L^{2}\left(\mathbb{R}\right)$ of) the
mapping\begin{eqnarray}
h & \rightarrow & \left(Zh\right)\left(t,\nu\right)=b^{-1/2}\sum_{k=-\infty}^{\infty}h\left(\frac{t-k}{b}\right)e^{2\pi ik\nu},\quad t,\nu\in\mathbb{R}.\end{eqnarray}
We refer to \cite{zz93} and to \cite[Sec. 1.5]{jans98}, for more
details on the Zak transform and its role in Gabor analysis.

For $f,h\in L^{2}\left(\mathbb{R}\right)$ we set (when $t,\nu\in\mathbb{R}$)\begin{eqnarray}
\Phi^{f}\left(t,\nu\right) & = & p^{-1/2}\left(\left(Zf\right)\left(t-l\frac{p}{q},\nu+\frac{k}{p}\right)\right)_{k=0,\ldots,p-1,\, l=0,\ldots,q-1},\end{eqnarray}
and \begin{eqnarray}
A^{f,h}\left(t,\nu\right) & = & \left(A_{k,r}^{f,h}\left(t,\nu\right)\right)_{k,r=0,\ldots,p-1}=\Phi^{f}\left(t,\nu\right)\left(\Phi^{h}\left(t,\nu\right)\right)^{*},\end{eqnarray}
where the $^{*}$ denotes conjugate transpose. Now $\left(g,a,b\right)$
is a Gabor frame, with frame bounds $A>0$, $B<\infty$ if and only
if we have $AI_{p\times p}\leq A^{gg}\left(t,\nu\right)\leq BI_{p\times p}$
for almost all $t,\nu\in\mathbb{R}$, with $A$ and $B$ the largest
and smallest positive real number for which the respective inequalities
hold. The frame operator $S$ of $\left(g,a,b\right)$ is {}``represented''
by $A^{gg}$ through the formula\begin{eqnarray}
\Phi^{Sf} & = & A^{gg}\Phi^{f},\quad f\in L^{2}\left(\mathbb{R}\right),\label{eq:ZZ-rep}\end{eqnarray}
 with matrix multiplication at each point $\left(t,\nu\right)\in\mathbb{R}$
on the right-hand side of (\ref{eq:ZZ-rep}). This formula extends
as follows. Assume that $\varphi$ is continuous and positive on $\left[A,B\right]$.
Then\begin{eqnarray}
\Phi^{\varphi\left(S\right)f} & = & \varphi\left(A^{gg}\right)\Phi^{f},\quad f\in L^{2}\left(\mathbb{R}\right),\label{eq:phi-ZZ}\end{eqnarray}
which is the basic formula for functional calculus in the Zak transform
domain.

We consider in this section the critical case $a=b=1$ (in Sections
\ref{sec:Implementational-aspects} and \ref{sec:Experiments} more
general rational $ab$ will be dealt with). Then considerable simplifications
occur since all the matrices $\Phi$, $A$ reduce to scalars. The
formula (\ref{eq:phi-ZZ}) then becomes\begin{eqnarray}
\left(Z\left(\varphi\left(S\right)f\right)\right)\left(t,\nu\right) & = & \varphi\left(\left|\left(Zg\right)\left(t,\nu\right)\right|^{2}\right)\left(Zf\right)\left(t,\nu\right),\quad t,\nu\in\mathbb{R},\end{eqnarray}
for $f\in L^{2}\left(\mathbb{R}\right)$. In particular we have\begin{eqnarray}
\left(Zg^{t}\right)\left(t,\nu\right) & = & \left(Z\left(S^{-1/2}g\right)\right)\left(t,\nu\right)=\frac{\left(Zg\right)\left(t,\nu\right)}{\left|\left(Zg\right)\left(t,\nu\right)\right|},\quad t,\nu\in\mathbb{R},\\
\left(Zg^{d}\right)\left(t,\nu\right) & = & \left(Z\left(S^{-1}g\right)\right)\left(t,\nu\right)=\frac{\left(Zg\right)\left(t,\nu\right)}{\left|\left(Zg\right)\left(t,\nu\right)\right|^{2}}\nonumber \\
 & = & \frac{1}{\left(Zg\right)^{*}\left(t,\nu\right)},\quad t,\nu\in\mathbb{R}.\end{eqnarray}

To illustrate the relevance for the algorithms, we consider algorithms
II and IV for all scaling strategies. As to initial scaling, we assume
that $g$ and $S$ are scaled such that $\left(g,a=1,b=1\right)$
has best upper frame bound $B<2$, which means that $\left|Zg\right|^{2}<2$
everywhere. We let\begin{eqnarray}
G=Zg & , & \Gamma_{k}=Z\gamma_{k}.\end{eqnarray}
Then by functional calculus in the Zak transform domain, the algorithms
II and IV (initial scaling) assume the form\begin{eqnarray}
\Gamma_{0}=G & ; & \Gamma_{k+1}=\frac{3}{2}\Gamma_{k}-\frac{1}{2}\left|\Gamma_{k}\right|^{2}\Gamma_{k},\quad k=0,1,\ldots,\label{eq:II-zak}\end{eqnarray}
and\begin{eqnarray}
\Gamma_{0}=G & ; & \Gamma_{k+1}=2\Gamma_{k}-\left|\Gamma_{k}\right|^{2}G,\quad k=0,1,\ldots,\label{eq:IV-zak}\end{eqnarray}
respectively, where the relations in (\ref{eq:II-zak}) and (\ref{eq:IV-zak})
are to be considered at each point $\left(t,\nu\right)\in\mathbb{R}$.
These recursions are then quite easily analyzed by elementary means.
For instance, one sees that the assumption $B<3$ is necessary and
sufficient for (\ref{eq:II-zak}) to converge to $\exp\left(i\arg\left(G\right)\right)$
everywhere, while the assumption $B<2$ is necessary and sufficient
for (\ref{eq:IV-zak}) to converge to $1/G^{*}$ everywhere. Unbounded
recursions result when we would have allowed $B$ to be larger than
$5$ and $2$, respectively.

Next we consider algorithms II, IV using norm scaling so that (\ref{eq:II-zak})
and (\ref{eq:IV-zak}) are to be replaced by\begin{eqnarray}
\Gamma_{0}=G & ; & \Gamma_{k+1}=\frac{3}{2}\frac{\Gamma_{k}}{\left\Vert \Gamma_{k}\right\Vert }-\frac{1}{2}\frac{\left|\Gamma_{k}\right|^{2}\Gamma_{k}}{\left\Vert \Gamma_{k}^{3}\right\Vert },\quad k=0,1,\ldots,\label{eq:II-zak-smart}\end{eqnarray}
and\begin{eqnarray}
\Gamma_{0}=G & ; & \Gamma_{k+1}=2\frac{\Gamma_{k}}{\left\Vert \Gamma_{k}\right\Vert }-\frac{\left|\Gamma_{k}\right|^{2}G}{\left\Vert \Gamma_{k}^{2}G\right\Vert },\quad k=0,1,\ldots\,.\label{eq:IV-zak-smart}\end{eqnarray}
The norms used here are $L^{2}\left(\left[0,1\right)^{2}\right)$-norms.
We consider the case that\begin{eqnarray}
Zg=1\textrm{ on }N & , & Zg=x>0\textrm{ on }M,\end{eqnarray}
where $N,M$ are two measurable sets $\subset[0,1)^{2}$ such that
$N\cap M=\emptyset$, $N\cup M=[0,1)^{2}$. Then $\left(g,a=1,b=1\right)$
is a Gabor frame with best frame bounds $A=\min\left(1,x^{2}\right)$,
$B=\max\left(1,x^{2}\right)$. Furthermore, \begin{eqnarray}
Zg^{t}=1\textrm{ on }[0,1)^{2} & ; & Zg^{d}=\frac{1}{x}\textrm{ on }M.\end{eqnarray}
We have for both algorithms II and IV that\begin{eqnarray}
\Gamma_{k}=c_{k}\textrm{ on }N & , & \Gamma_{k}=d_{k}\textrm{ on }M,\label{eq:cd-on-NM}\end{eqnarray}
where $c_{k}$, $d_{k}$ follow recursions that can be made completely
explicit (using that, for instance, $\left\Vert \Gamma_{k}\right\Vert =\left(\left(1-\varepsilon\right)c_{k}^{2}+\varepsilon d_{k}^{2}\right)^{1/2}$,
where $\varepsilon=\mu\left(M\right)$). Due to the norming operations
in the recursion steps, either recursion stays bounded.

We consider the case that $\varepsilon=\mu\left(M\right)\ll1$. Then
an elementary analysis shows the following: There is a $\delta>0$
such that for recursion (\ref{eq:II-zak-smart}), (\ref{eq:cd-on-NM})
there holds

\begin{itemize}
\item $x\in\left(0,\sqrt{3}-\delta\right)\Rightarrow c_{k},d_{k}\rightarrow1$,
\item $x\in\left(\sqrt{3}+\delta,\sqrt{5}-\delta\right)\Rightarrow c_{k}\rightarrow1,d_{k}\rightarrow-1$
\item $x\in\left(\sqrt{5}+\delta,\infty\right)\Rightarrow\textrm{ chaotic behaviour}$.
\end{itemize}
Also, there is a $\delta>0$ such that for recursion (\ref{eq:IV-zak-smart}),
(\ref{eq:cd-on-NM}) there holds

\begin{itemize}
\item $x\in\left(0,\sqrt{2}-\delta\right)\Rightarrow c_{k}\rightarrow1,d_{k}\rightarrow\frac{1}{x}$,
\item $x\in\left(\sqrt{2}+\delta,\infty\right)\Rightarrow d_{k}<0$.
\end{itemize}

\section{\label{sec:Sampling-and-periodization}Sampling and periodization
of Gabor frames}

The algorithms considered in this paper and in \cite{jans02,jans02b,jans03}
have been formulated for time-continuous Gabor frames while the tests
we perform must take place in a finite setting. The transition from
continuous to discrete/finite Gabor frames by sampling and periodization
has been discussed in \cite{Jans97} and later in \cite{kai05,sonder06},
see also \cite[Subsec. 8.4]{chrstr03} and \cite[Secs. 10.2 and 10.3]{Chr03}.
Let $a$, $M$ be positive integers, and assume that $\left(g,a,1/M\right)$
is a Gabor frame with frame bounds $A>0$, $B<\infty$. Furthermore,
assume that $g$ satisfies the aforementioned condition A and the
so-called condition R:

\begin{lyxlist}{MMM.}
\item [R:]\begin{eqnarray}
\lim_{\varepsilon\rightarrow0}\sum_{j=-\infty}^{\infty}\frac{1}{\varepsilon}\int_{-\varepsilon/2}^{\varepsilon/2}\left|g\left(j+u\right)-g\left(j\right)\right|^{2}du & = & 0.\end{eqnarray}

\end{lyxlist}
The conditions R and A are not very restrictive; they are, for instance,
satisfied by all members $g$ of Feichtinger's algebra $S_{0}$, see
\cite[comment after Thm. 8.4.2]{chrstr03}. Then the system\begin{eqnarray}
g_{na,m/M}^{D} & := & \left(g_{na,m/M}\left(r\right)\right)_{r\in\mathbb{Z}},\quad n\in\mathbb{Z},m=0,\ldots,M-1\label{eq:discr-gab-syst}\end{eqnarray}
 is a discrete Gabor frame with frame bounds $A$, $B$, the dual
window $g^{d}$ of the frame $\left(g,a,1/M\right)$ satisfies conditions
R and A, and the dual window $\left(g^{D}\right)^{d}$ corresponding
to the discrete Gabor system in (\ref{eq:discr-gab-syst}) is obtained
by sampling $g^{d}$:\begin{eqnarray}
\left(g^{D}\right)^{d}\left(r\right) & = & \left(g^{d}\right)^{D}\left(r\right),\quad r\in\mathbb{Z}.\end{eqnarray}

The transition from discrete Gabor frames to discrete, periodic Gabor
frames is just as convenient. Assume that we have a $g\in l^{1}\left(\mathbb{Z}\right)$
such that the discrete Gabor system $\left(g_{na,m/M}\right)_{n\in\mathbb{Z},m=0,\ldots,M-1}$
is a discrete Gabor frame with frame bounds $A>0$ $B<\infty$. Let
$L=Na=Mb$ for some positive integers $N,b$, and define\begin{eqnarray}
g^{P}\left(r\right) & = & \sum_{j=-\infty}^{\infty}g\left(r-jL\right),\quad r\in\mathbb{Z}.\end{eqnarray}
Then the system\begin{equation}
\left(g_{an,m/M}^{P}\left(r\right)\right)_{r\in\mathbb{Z}},\quad n=0,\ldots,N-1,m=0,\ldots,M-1,\label{eq:disc-per-gab-sys}\end{equation}
 is a discrete, periodic Gabor system with frame bounds $A$, $B$,
the dual window $g^{d}$ of the discrete Gabor system is in $l^{1}\left(\mathbb{Z}\right)$,
and the dual window $\left(g^{P}\right)^{d}$ corresponding to the
discrete periodic Gabor system in (\ref{eq:disc-per-gab-sys}) is
obtained by periodizing $g^{d}$:\begin{eqnarray}
\left(g^{P}\right)^{d}\left(r\right) & = & \left(g^{d}\right)^{P}\left(r\right)=\sum_{j=-\infty}^{\infty}g^{d}\left(r-jL\right),\quad r\in\mathbb{Z}.\end{eqnarray}
An important extension of these results is given in \cite[Subsec. 8.4]{chrstr03}.
Assume that $\varphi$ is analytic in an open neighbourhood containing
$\left[A,B\right]$, where $A>0$, $B<\infty$ are frame bounds of
the Gabor frame $\left(g,a,1/M\right)$ with $g$ satisfying condition
R and A. Then $\varphi\left(S\right)g$ satisfies R and A as well,
and\begin{eqnarray}
\left(\varphi\left(S\right)g\right)^{D} & = & \varphi\left(S^{D}\right)g^{D},\end{eqnarray}
where $S^{D}$ is the frame operator corresponding to the system in
(\ref{eq:discr-gab-syst}). The approach in \cite[Subsec. 8.4]{chrstr03}
(which uses the Dunford representation of operators as well as theorems
of the Wiener $1/f$-type) can be mimicked so as to generalize the
transition result from discrete Gabor systems as above with $g\in l^{1}\left(\mathbb{Z}\right)$
to discrete, periodic Gabor systems. Thus, with $\varphi$ as above
and $\left(g_{na,m/M}\right)_{n\in\mathbb{Z},m=0,\ldots,M-1}$ a discrete
Gabor system with $g\in l^{1}\left(\mathbb{Z}\right)$, frame bounds
$A>0$, $B<\infty$ and frame operator $S$, we have $\varphi\left(S\right)g\in l^{1}\left(\mathbb{Z}\right)$
and \begin{eqnarray}
\left(\varphi\left(S\right)g\right)^{P} & = & \varphi\left(S^{P}\right)g^{P},\end{eqnarray}
where $S^{P}$ is the frame operator of the system in (\ref{eq:disc-per-gab-sys}).
In particular, we see that the sampling-and-periodization approach
is valid for the tight window $g^{t}$ in which case we should consider
$\varphi\left(s\right)=s^{-1/2}$.

It follows from the above results that the algorithms can be considered
for discrete and for discrete, periodic Gabor frames. The findings
for these systems are of direct relevance to the algorithms we have
considered for the time-continuous case.

\section{\label{sec:Implementational-aspects}Implementational aspects}

All implementations are done in the finite, discrete setting of Gabor
frames. We denote for $g\in\mathbb{C}^{L}$ and $a,b\in\mathbb{N}$
by $\left(g,a,b\right)$ the collection of time-frequency shifted
windows\begin{equation}
g_{na,mb},\quad n\in\mathbb{Z},m\in\mathbb{Z},\end{equation}
 where for $j,k\in\mathbb{Z}$ we denote\begin{equation}
g_{j,k}=e^{2\pi ikl/L}g(l-j),\quad l=0,\ldots L-1.\end{equation}
Note that it must hold that $L=Na=Mb$ for some $M,N\in\mathbb{N}.$
Additionally, we define $c,d,p,q\in\mathbb{N}$ by \begin{equation}
c=\gcd\left(a,M\right),\quad d=\gcd\left(b,N\right),\quad p=\frac{a}{c}=\frac{b}{d},\quad q=\frac{M}{c}=\frac{N}{d}.\end{equation}
With these numbers, the density of the Gabor system can be written
as$\left(ab\right)/L=p/q,$ where $p/q$ is a irreducible fraction.
It holds that $L=cdpq$.

\subsection{Matrix representation and the SVD}

Let $O_{g}\in\mathbb{C}^{L\times MN}$ be the matrix representation
of the synthesis operator of a Gabor frame so that\begin{eqnarray}
\left(O_{g}\right)_{l,m+nM} & = & g_{ma,nb}\left(l\right),\quad l=0,\ldots,L-1,\end{eqnarray}
for $m=0,\ldots,M-1$, $n=0,\ldots,N-1$. Hence $O_{g}$ has the column
vectors $g_{ma,nb}$. The matrix representation of the frame operator
corresponding to $\left(g,a,b\right)$ is then given as $O_{g}O_{g}^{*}$.
Since $\left(g,a,b\right)$ is a frame we have that $O_{g}$ has full
rank $L\leq MN$.

Assume that $\varphi$ is continuous an positive on $\sigma\left(S\right)$.
From\begin{eqnarray}
\left(\varphi\left(S\right)g\right)_{na,mb} & = & \varphi\left(S\right)g_{na,mb},\end{eqnarray}
 we have that\begin{eqnarray}
O_{\varphi\left(S\right)g} & = & \varphi\left(S\right)O_{g}=\varphi\left(O_{g}O_{g}^{*}\right)O_{g}.\end{eqnarray}
Furthermore, note that for the Frobenius norm $\left\Vert O_{g}\right\Vert _{\text{fro}}$
we have $\left\Vert O_{g}\right\Vert _{\text{fro}}^{2}=MN\left\Vert g\right\Vert ^{2},$
since all columns of $O_{g}$ have norm $\left\Vert g\right\Vert $.

The iterations can be written in terms of the synthesis operator matrices
as follows. Denote the synthesis operator matrix corresponding to
the Gabor frame $\left(\gamma_{k},a,b\right)$ by $\Omega_{k}$. Then
we can write the iteration step for algorithm II with norm scaling
as\begin{equation}
\Omega_{0}=O_{g};\,\Omega_{k+1}=\frac{3}{2}\frac{\Omega_{k}}{\left\Vert \Omega_{k}\right\Vert _{\text{fro}}}-\frac{1}{2}\frac{\left(\Omega_{k}\Omega^{*}\right)\Omega_{k}}{\left\Vert \left(\Omega_{k}\Omega_{k}^{*}\right)\Omega_{k}\right\Vert _{\text{fro}}},\quad k=0,1,\ldots.\label{eq:II-by-gabmat}\end{equation}

We shall consider the \emph{thin} SVD of the synthesis operator matrices.
Thus we let $O_{g}=U\Sigma V^{*}$, where $U\in\mathbb{C}^{L\times L}$
is unitary ($O_{g}$ has full rank), $\Sigma\in\mathbb{R}^{L\times L}$
is a diagonal matrix with positive diagonal elements and $V\in\mathbb{C}^{MN\times L}$
has orthonormal columns. With $\varphi$ as above, we compute the
thin SVD of $O_{\varphi\left(S\right)g}$ as 

\begin{eqnarray}
O_{\varphi\left(S\right)g} & = & \varphi\left(O_{g}O_{g}^{*}\right)O_{g}=\varphi\left(U\Sigma^{2}U^{*}\right)U\Sigma V^{*}\nonumber \\
 & = & U\varphi\left(\Sigma^{2}\right)U^{*}U\Sigma V^{*}=U\Sigma\varphi\left(\Sigma^{2}\right)V^{*}.\label{eq:phi-SVD}\end{eqnarray}
Here we have used that $\varphi\left(U\Sigma^{2}U^{*}\right)=U\varphi\left(\Sigma^{2}\right)U^{*}$,
a basic fact in the functional calculus of matrices. The equation
(\ref{eq:phi-SVD}) shows that $O_{\varphi\left(S\right)g}$ has the
same right and left singular vectors as $O_{g}$, and the singular
values transform according to $\varphi\rightarrow\sigma\varphi^{2}\left(\sigma\right)$.
As a consequence we have, \begin{eqnarray}
O_{g^{t}}=UV^{*} & , & O_{g^{d}}=U\Sigma^{-1}V^{*},\label{eq:UV-SVD}\end{eqnarray}
for the synthesis operators corresponding to $\left(g^{t},a,b\right)$
and $\left(g^{d},a,b\right)$, respectively, for which we should take
$\varphi\left(s\right)=s^{-1/2}$ and $\varphi\left(s\right)=s^{-1}$
in (\ref{eq:phi-SVD}). We thus see that we have obtained the matrices
occurring in the polar decomposition of $O_{g}$ and the Moore-Penrose
pseudo-inverse of $O_{g}$.

A further observation is that $\left\Vert O_{g}\right\Vert _{\text{fro}}^{2}=\sum_{j=1}^{j=L}\sigma_{j}^{2},$
where $\sigma_{j}$, $j=0,\ldots,L-1$, are the singular values of
$O_{g}$. Letting $\varphi_{k,j}$ be the singular values of $\Omega_{k}$
and using that\[
\Omega_{k}=U\Sigma_{k}V^{*},\]
we can write the iteration step in (\ref{eq:II-by-gabmat}) on the
level of singular values as

\begin{eqnarray}
\sigma_{k+1,j} & = & \frac{3}{2}\alpha_{k}\sigma_{k,j}-\frac{1}{2}\beta_{k}\sigma_{k,j}^{3},\quad\alpha_{k}=\frac{1}{\sqrt{\sum_{j=1}^{j=L}\sigma_{k,j}^{2}}},\,\beta_{k}=\frac{1}{\sqrt{\sum_{j=1}^{j=L}\sigma_{k,j}^{6}}},\label{eq:scal-eq-norm}\end{eqnarray}
where $j=0,\ldots,L-1,$ and $k=0,1,\ldots$.

\subsection{Factorization of finite, discrete Gabor systems}

Similar to the Zibulski-Zeevi representation of the Gabor frame operator
in the continuous case, see \cite{zz93,zz97,jans98}, it is possible
to compute the actions of the finite, discrete Gabor frame operator
(and also the analysis and synthesis operators) very efficiently.
Several equivalent methods exists using almost the same number of
operations, but differing in the order. In \cite{bastgeil96} a finite,
discrete version of the Zibulski-Zeevi representation is developed.
Another method was developed in \cite{prinz96} and \cite{Str98a}.
Unfortunately, \cite{Str98a} contains some errors, which have been
corrected in \cite{bal05}. In the following we shall present the
Zak-transform method from \cite{bastgeil96}.

For $h\in\mathbb{C}^{L}$ and $K\in\left\{ 0,...,L-1\right\} $ such
that $\frac{L}{K}\in\mathbb{N}$, we define the finite, discrete Zak
transform $Z_{K}h$ by \begin{eqnarray}
\left(Z_{K}h\right)\left(r,s\right) & = & \sqrt{\frac{K}{L}}\sum_{l=0}^{L/K-1}h\left(r-lK\right)e^{2\pi islK/L},\quad r,s\in\mathbb{Z}.\end{eqnarray}
The finite, discrete Zak transform is quasi-periodic in its first
variable and periodic in the second,\begin{eqnarray}
\left(Z_{K}h\right)\left(r+kK,s+l\frac{K}{L}\right) & = & e^{2\pi iksK/L}\left(Z_{K}h\right)\left(r,s\right),\label{eq:FDZAK-quasi}\end{eqnarray}
see \cite{jans94b} for more details. The values $\left(Z_{K}h\right)\left(r,s\right)$
of a finite, discrete Zak-transform on the fundamental domain $r=0,\ldots,K-1$,
$s=0,\ldots,K/L-1$ can be calculated efficiently by $K$ FFT's of
length $K/L$. To obtain values outside the fundamental domain, the
quasi-periodicity relation (\ref{eq:FDZAK-quasi}) can be used.

We define the $cd$ matrices $\Phi_{r,s}^{f}$ of size $p\times q$
and the $p\times p$ matrices $A_{r,s}^{f}$ by \begin{eqnarray}
\Phi_{r,s}^{f} & =\sqrt{cdq}\left(\left(Z_{a}f\right)\left(r+kM,s+ld\right)\right)_{k=0,...,p-1;\, l=0,...,q-1} & ,\label{eq:fac}\end{eqnarray}
where $r=0,\ldots,c-1$, $s=0,\ldots,d-1$ and \begin{eqnarray}
A_{r,s}^{f,h} & = & \left(A_{r,s}^{f,h}\right)_{k,l=0,\ldots,p-1}=\Phi_{r,s}^{f}\left(\Phi_{r,s}^{h}\right)^{*}.\end{eqnarray}
With these definitions is holds that the frame operator $S$ of $\left(g,a,b\right)$
is {}``represented'' by $A^{gg}$ through the formula\begin{eqnarray}
\Phi^{Sf} & = & A^{gg}\Phi^{f},\quad f\in\mathbb{C}^{L},\label{eq:ZZ-rep-D}\end{eqnarray}
see \cite{bastgeil96}.

With this efficient representation of the frame operator of a finite,
discrete Gabor system, we may express the iterations schemes in the
finite, discrete Zak domain. We let\begin{eqnarray}
G=\Phi^{g} & , & \Gamma_{k}=\Phi^{\gamma_{k}},\, A_{k}=\Phi^{\gamma_{k},\gamma_{k}}=\Gamma_{k}\left(\Gamma_{k}\right)^{*}.\end{eqnarray}
By functional calculus, algorithm II in the finite, discrete Zak transform
takes the form\begin{eqnarray}
\Gamma_{0}=G & ; & \Gamma_{k+1}=\frac{3}{2}\Gamma_{k}-\frac{1}{2}A_{k}\Gamma_{k},\quad k=0,1,\ldots.\label{eq:II-zak-D}\end{eqnarray}
The expressions for the other iterations types are similar.

\subsection{Other methods}

We have considered two other methods of computing the canonical tight
window utilizing the factorization (\ref{eq:fac}). To calculate the
factorization of the canonical tight window $g^{t}$,$\Phi^{g^{t}}$,
we use an eigenvalue decomposition of the factorization of the frame
operator of $\left(g,a,b\right)$: For each $r=0,\ldots,c-1$, $s=0,\ldots,d-1$
compute $U_{r,s},D_{r,s}$ such that $A_{r,s}^{gg}=U_{r,s}D_{r,s}U_{r,s}^{*},$
where $U_{r,s}$ is unitary and $D_{r,s}$ is diagonal and set\begin{eqnarray}
\Phi_{r,s}^{g^{t}} & = & U_{r,s}D_{r,s}^{-1/2}U_{r,s}^{*}\Phi_{r,s}^{g}.\end{eqnarray}
 We shall refer to this method as the EIG method. The other method
uses (\ref{eq:UV-SVD}) applied to the matrices of the factorization:
For each $r=0,\ldots,c-1$, $s=0,\ldots,d-1$ compute $U_{r,s},D_{r,s},V_{r,s}$
such that $\Phi_{r,s}^{g}=U_{r,s}D_{r,s}V_{r,s}^{*},$ where $U_{r,s}$
is unitary, $D_{r,s}$ is diagonal and $V_{r,s}$ has orthonormal
columns. Then it follows from functional calculus in the Zak transform
domain (pretty much as in (\ref{eq:phi-SVD}); also see (\ref{eq:phi-ZZ}))
that \begin{eqnarray}
\Phi_{r,s}^{g^{t}} & = & U_{r,s}V_{r,s}^{*}.\end{eqnarray}
We shall refer to this method as the SVD method.

For computing the canonical dual window we have considered simply
inverting the matrices of the factorization of the frame operator:\begin{eqnarray*}
\Phi_{r,s}^{g^{d}} & = & \left(A_{r,s}^{gg}\right)^{-1}\Phi_{r,s}^{g}.\end{eqnarray*}
We shall refer to this as the INV method.

\subsection{Implementational costs}

The computation of $\Phi^{g}$ needs to be done before the iteration
step. It can be computed using $5L\log_{2}d$ flops. This transforms
the initial window $g$ into the finite, discrete Zak domain. All
computations in this domain are then done by multiplication of $p\times q$
and $p\times p$ matrices. The transform $\Phi$ is unitary from $\mathbb{C}^{L}$
with Euclidean norm into $\mathbb{C}^{c\times d\times p\times q}$,
also with Euclidean norm. This gives an easy way to calculate the
norms needed for the norm scaling.

We count the number of real floating point operation needed, and assume
that everything is done using complex arithmetics. The flop count
for a single iteration step in the transform domain for each of the
5 algorithms can be seen in Table \ref{cap:flopcount}.

\begin{table}

\caption{\label{cap:flopcount}}

\begin{longtable}{p{2cm}|p{6cm}}
\hline 
Method:&
Flop count per iteration:\tabularnewline
\hline
I.&
$16Lp+\frac{4}{3}cdp^{3}$.\tabularnewline
II.&
$16Lp$. \tabularnewline
III.&
$24Lp$.\tabularnewline
IV.&
$16Lp$.\tabularnewline
V.&
$24Lp+8cdp^{3}$.\tabularnewline
\hline 
&
Total flop count:\tabularnewline
\hline 
INV.&
$16Lp+\frac{4}{3}cdp^{3}$.\tabularnewline
EIG.&
$24Lp+14cdp^{3}$.\tabularnewline
SVD.&
$64Lp+32cdp^{3}$.\tabularnewline
\hline
\end{longtable}

This table shows the flop count of each of the considered methods.
The flop count does not include the cost of the pre- and post-factorization.
The application of an inverse matrix needed for the algorithms I and
INV is done using a Cholesky factorization followed by two substitutions.
An iteration step of V takes more flops to compute than an iteration
step of III, because we need to compute the two terms $S_{k}g$ and
$S_{k}\gamma_{k}$. The flop counts for EIG and SVD methods are only
approximations, because eigenvalues and singular values can be calculated
by many different methods with different flop counts, and because
the process usually involves an iterative step, see \cite{golloan96}.\\

\end{table}

A quick comparison show that the iterative methods for computing the
tight window are comparable in number of flops to the EIG and SVD
methods, if the number of necessary iterations is not to big. For
the inverse iterations, the situation is different: Computing the
inverse of the block matrices by a direct approach requires only slightly
more flops than a single iteration step of algorithm IV, so an iterative
method will always use more flops than the direct approach. However,
there might be situations were it is not desirable to compute the
inverse. For instance, if the initial window $g$ has small support
then the iteration steps can be performed by multiple passes through
a filter bank.

\subsection{Stopping criterion}

Because of the guaranteed quadratic/cubic convergence of the algorithms,
it is possible to devise a simple yet powerful stopping criterion:
We consider the difference \begin{equation}
\frac{\left\Vert \gamma_{k+1}-\gamma_{k}\right\Vert }{\left\Vert \gamma_{k+1}\right\Vert }.\label{eq:rel-diff}\end{equation}
When this difference is close to the machine precision $eps$, the
algorithm considered has converged. This is a standard stopping criteria,
but using it this way means that we have done exactly one iteration
step too much. Therefore, we stop when (\ref{eq:rel-diff}) is less
that $\sqrt{eps}$ and $\sqrt[3]{eps}$ for the algorithms having
quadratic and cubic convergence, respectively.

\subsection{Window functions}

\begin{figure}
\includegraphics[%
  width=0.5\textwidth,
  keepaspectratio]{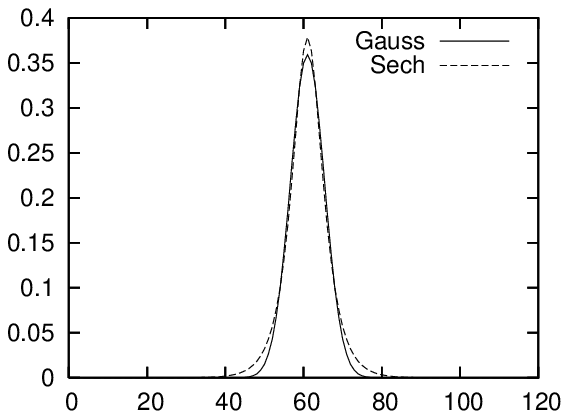}\includegraphics[%
  width=0.5\textwidth,
  keepaspectratio]{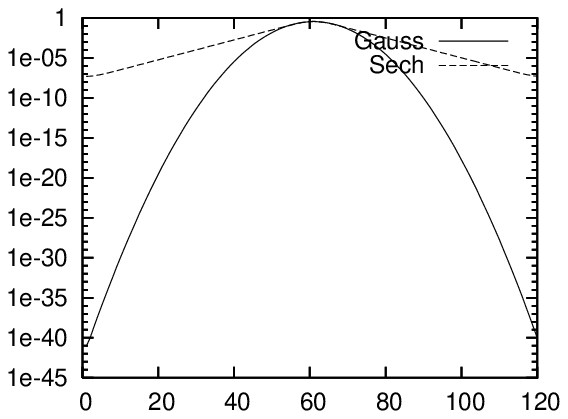}\\
\includegraphics[%
  width=0.5\textwidth,
  keepaspectratio]{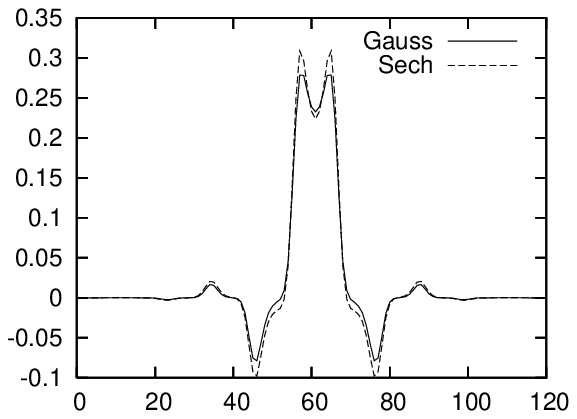}\includegraphics[%
  width=0.5\textwidth,
  keepaspectratio]{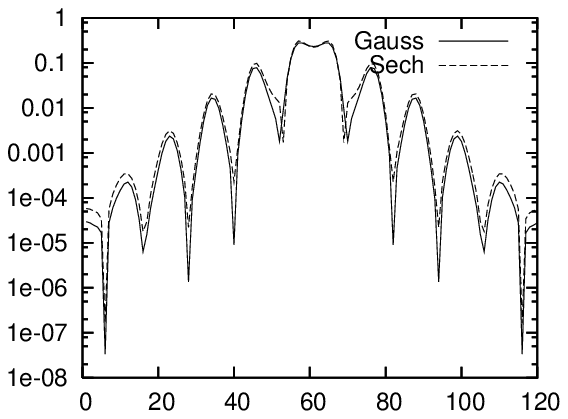}\\
\includegraphics[%
  width=0.5\textwidth,
  keepaspectratio]{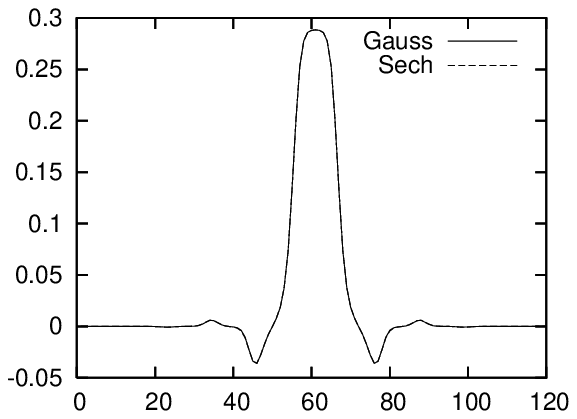}\includegraphics[%
  width=0.5\textwidth,
  keepaspectratio]{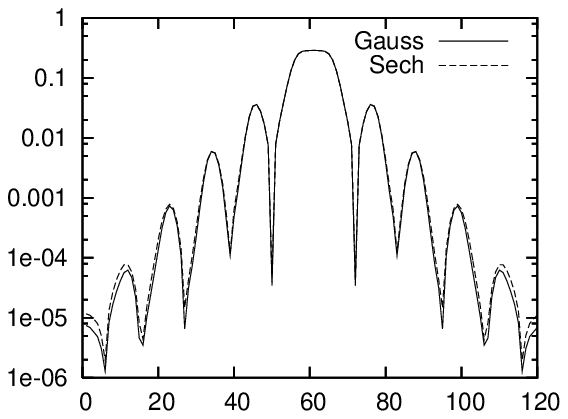}

\caption{\label{cap:windows}The first row of graphs shows the Gaussian function
$\varphi_{1}^{D}$ and the hyperbolic secant $\psi_{1}^{D}$. Both
functions have length $L=120$ and unit norm. The second row of graphs
shows their canonical dual windows, and the third row shows their
canonical tight windows. In the second column, the functions are shown
on a logarithmic scale to display their decay.\protect \\
}
\end{figure}

As the basic window functions we shall use the Gaussian $\varphi\in L^{2}\left(\mathbb{R}\right)$
and the hyperbolic secant $\psi\in L^{2}\left(\mathbb{R}\right)$
given by\begin{eqnarray}
\varphi\left(t\right) & = & \left(\frac{1}{2}\right)^{-1/4}e^{-\pi t^{2}},\quad t\in\mathbb{R},\\
\psi\left(t\right) & = & \sqrt{\frac{\pi}{2}}\textrm{sech}\left(\pi t\right),\quad t\in\mathbb{R}.\end{eqnarray}
Both functions are invariant with respect to Fourier transformation.
In order to generate a range of functions, we introduce a parameter
$w$ that dilates the functions by the unitary operator $D_{w}:\, L^{2}\left(\mathbb{R}\right)\rightarrow L^{2}\left(\mathbb{R}\right)$
given by \begin{eqnarray}
\left(D_{w}f\right)\left(t\right) & =\left(w\right)^{-1/4} & f\left(\frac{t}{\sqrt{w}}\right),\quad t\in\mathbb{R}.\end{eqnarray}
Applying this gives\begin{eqnarray}
\left(D_{w}\varphi\right)\left(t\right) & = & \varphi_{w}\left(t\right)=\left(\frac{w}{2}\right)^{-1/4}e^{-\pi t^{2}/w},\quad t\in\mathbb{R},\\
\left(D_{w}\psi\right)\left(t\right) & = & \psi_{w}\left(t\right)=\sqrt{\frac{\pi}{2}}w^{-1/4}\textrm{sech}\left(t\frac{\pi}{\sqrt{w}}\right),\quad t\in\mathbb{R}.\end{eqnarray}
It holds that the Fourier transform of $\varphi_{w}$ is $\varphi_{1/w}$
and similarly for $\psi_{w}$. As window functions for the testing
of the iterative algorithms we shall use finite, discrete versions
of these, obtained by the sampling-and-periodization process described
in Sec. \ref{sec:Sampling-and-periodization}:\begin{eqnarray}
\varphi_{w}^{D}(l) & = & \left(\frac{wL}{2}\right)^{-1/4}\sum_{k\in\mathbb{Z}}e^{-\pi\left(l/\sqrt{L}-k\sqrt{L}\right)^{2}/w},\quad l=0,\ldots,L-1,\\
\psi_{w}^{D}\left(l\right) & = & \sqrt{\frac{\pi}{2}}\left(wL\right)^{-1/4}\left(\frac{2\sqrt{wL}}{\pi}\right)^{-1/2}\sum_{k\in\mathbb{Z}}\textrm{sech}\left(\left(\frac{l}{\sqrt{L}}-k\sqrt{L}\right)\frac{\pi}{\sqrt{w}}\right).\end{eqnarray}
The properties from the continuous setting carry over: The functions
have unit norm, and the Discrete Fourier Transform of $\varphi_{w}^{D}$
is $\varphi_{1/w}^{D}$ and similarly for $\psi_{w}^{D}$. For more
details on the hyperbolic secant as a Gabor window, see \cite{jansstro02}.
Figure \ref{cap:windows} shows a Gaussian and a hyperbolic secant
and their respective canonical dual and tight windows. The different
decay properties of the two bell-shaped functions are visible on a
logarithmic scale. However, even though the Gaussian and the hyperbolic
secant have different decay properties, their canonical windows have
almost the same.

\begin{figure}
\subfigure[MONSTER]{\includegraphics[%
  width=0.5\textwidth,
  keepaspectratio]{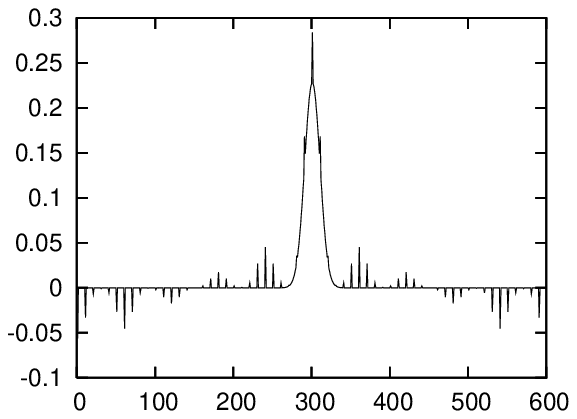}}\subfigure[DFT of MONSTER]{\includegraphics[%
  width=0.5\textwidth,
  keepaspectratio]{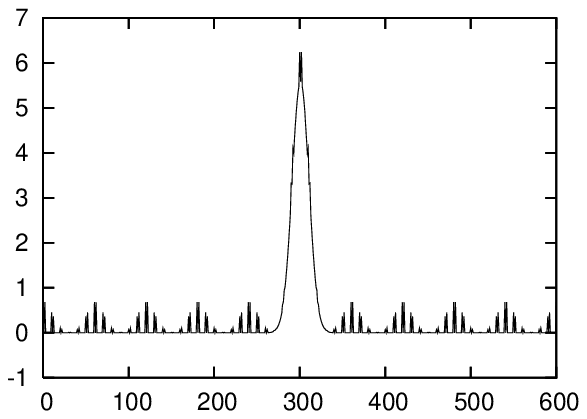}}

\caption{\label{cap:jmonster}Fig. (a) shows the function MONSTER of length
$600$ with a single singular value set to $\sigma_{\text{real}}=6$.
Fig. (b) shows the Discrete Fourier Transform of the function.\protect \\
}
\end{figure}
To produce examples for which the norm scaling methods diverges, we
have constructed a function (MONSTER) which is a Gaussian function
modified in such a way that the first singular value, $\sigma_{\text{real}}$,
of the matrix representation of the Gabor synthesis operator that
corresponds to a real and symmetric singular vector, is given a large
value. The function is shown on Figure \ref{cap:convergence}. This
function is a generalization to the case of rational oversampling
of the counterexample given in Sec. \ref{sec:Zak-domain-considerations}
and exploits that the iterations can be considered as scalar iterations
of the singular values of the Gabor synthesis operator, (\ref{eq:scal-eq-norm}).

\section{\label{sec:Experiments}Experiments}

This section contains the results from the experiments we have done
in order to test the algorithms thoroughly and to demonstrate the
various aspects of the algorithms that have been shown analytically.
The computations have been done in Matlab and Octave, and the full
source code is available for download from \url{http://www2.mat.dtu.dk/people/P.Soendergaard/iteralg/}.
We will show figures demonstrating the important aspects, but since
we cannot include all material, the reader is encouraged to download
the software and experiment.

\subsection{Convergence and divergence of norm scaling}

\begin{figure}
\subfigure[Tight iterations.]{\includegraphics[%
  width=0.5\textwidth,
  keepaspectratio]{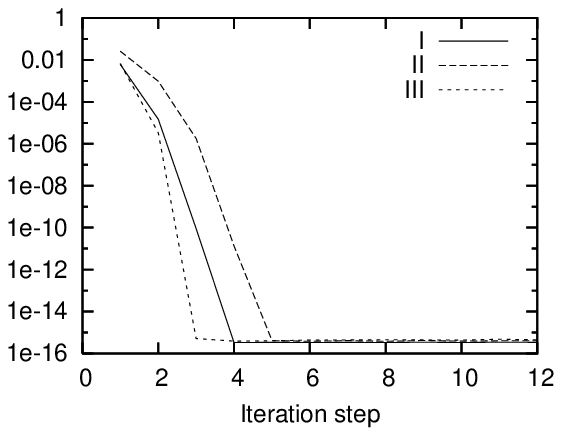}}\subfigure[Dual iterations.]{\includegraphics[%
  width=0.5\textwidth,
  keepaspectratio]{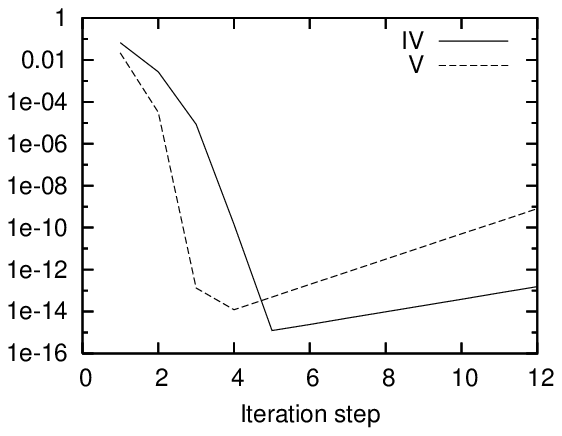}}

\caption{\label{cap:convergence}The figure shows the behaviour of the 5 iteration
types for the first 12 iterations of each. The y-axis shows the $l^{2}$-norm
of the difference between the iteration step and a precomputed, normalized
solution. The system considered in the Gabor frame for $\mathbb{C}^{432}$,
$\left(\varphi_{1}^{D},18,18\right)$. The Gabor frame has a frame
bound ratio of $B/A=2.03$.\protect \\
}
\end{figure}

Figure \ref{cap:convergence} show the convergence behaviour for a
well-conditioned problem. The figure shows that algorithm I,II and
IV exhibit quadratic convergence, III and V exhibit cubic convergence
as proved in Subsections \ref{sub:Iterations-for-t} and \ref{sub:Iterations-for-d}.
Furthermore, the algorithms for computing the tight window stay converged
close to the machine precision, while the algorithms for computing
the dual window diverge. Algorithm V also diverges faster than IV.
This is as proved in Subsection \ref{sub:Influence-of-out-of-space};
the slopes of the two line segments beyond the $5\mathrm{^{th}}$
iteration in Fig. \ref{cap:convergence}(b) corresponds to divergence
factors 2 (for IV) and 4 (for V). A visible numerical aspect is that
iteration V is not able to reach full precision, because the iterand
is quickly affected by the buildup of numerical errors. The convergence
behaviour of the algorithms for the initial window being a hyperbolic
secant is almost the same.

\begin{figure}
\subfigure[Gaussian]{\includegraphics[%
  width=0.5\textwidth,
  keepaspectratio]{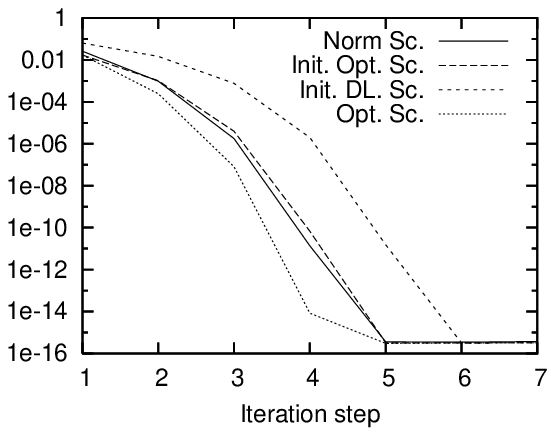}}\subfigure[Badly scaled Gaussian]{\includegraphics[%
  width=0.5\textwidth,
  keepaspectratio]{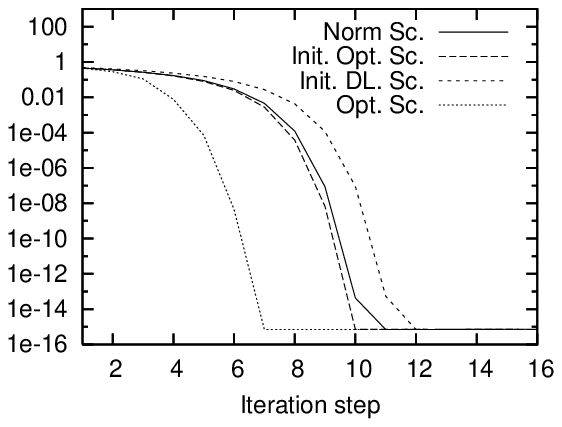}}

\caption{\label{cap:conv_scal}Fig. (a) shows the convergence behaviour for
algorithm II using four different scaling strategies: Norm scaling,
initial scaling by the optimal constant, initial scaling by the dual
lattice norm and constant optimal scaling. The system considered in
the Gabor frame for $\mathbb{C}^{432}$, $\left(\varphi_{1}^{D},18,18\right)$.
Fig (b) shows the same, but instead using the window $\varphi_{1/5}^{D}$.
This is a very narrow window, and the generated Gabor frame has a
frame bound ratio of $B/A=180.8$. \protect \\
}
\end{figure}

Two examples of using different scaling strategies is show on Figure
\ref{cap:conv_scal}(a) and \ref{cap:conv_scal}(b). The figures show
that initially scaling by the best scaling constant is as good as
using norm scaling and using initial scaling by an easily computable
scaling constant results in only 1-2 more iterations than using norm
scaling. Comparing these methods to the method using optimal scaling,
we see that for a well-conditioned problem then norm-scaling and optimal
initial scaling are close to the optimal convergence. For a worse
conditioned problem (Fig. \ref{cap:conv_scal}(b)), optimal scaling
clearly outperforms the other methods. However, this observation has
little practical relevance, because the computed canonical windows
$g^{d}$ and $g^{t}$ will have a bad time-frequency localization.
Higham \cite{high86p} uses a scaling strategy for algorithm I that
approximates the optimal scaling. This requires an estimate for the
smallest eigenvalue of the matrix, but this is easy to obtain since
the matrix is inverted as part of the iteration step. For algorithms
II-V we cannot use inversions, and so an estimate for the smallest
eigenvalue (or lower frame bound) is difficult to obtain. We have
therefore not pursued such a method for algorithms II-V.

\begin{figure}
\subfigure[II]{\includegraphics[%
  width=0.5\textwidth,
  keepaspectratio]{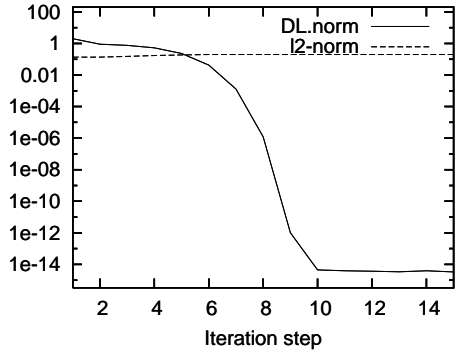}}\subfigure[IV]{\includegraphics[%
  width=0.5\textwidth,
  keepaspectratio]{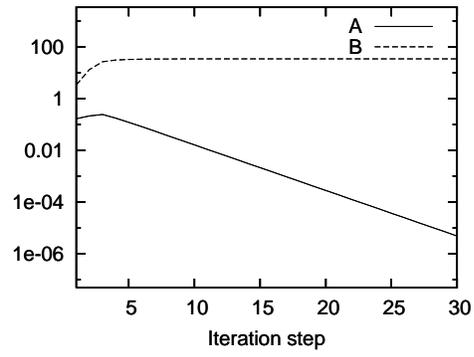}}

\caption{\label{cap:divergence}Fig. (a) shows the behaviour of the dual lattice
norm and the $l^{2}$-norm of the difference between the iteration
step and the normalized solution for a run of algorithm II using norm
scaling. Fig. (b) shows the behaviour of the best upper and lower
frame bound of $Z_{k}$ in each iteration step for a run of algorithm
IV using norm scaling. The system considered is in both cases are
the Gabor frame for $\mathbb{C}^{600}$ with $a=b=20$ using the MONSTER
function.\protect \\
}
\end{figure}

The iterations for computing the tight window are very robust when
using norm scaling. It is easy to create examples of Gabor systems
with frame bound ratios $B/A>10^{12}$ for which the iterations converge,
by using badly dilated Gaussians or by using a constant function with
a small amount of noise added. However, by using the $MONSTER$ function,
it is possible to create an example for which the norm scaling iterations
diverge. The behaviour of the dual lattice norm and $\left\Vert g-\gamma_{k}\right\Vert _{2}$
in each iteration step for a run of algorithm II step is shown in
Fig. \ref{cap:divergence}(a). It can be seen that the iteration converges
to the wrong tight window. Another typical behaviour is that the iteration
oscillates between two different functions with the same dual lattice
norm. The behaviour of algorithm IV on the same examples is shown
in Fig. \ref{cap:divergence}(b), the figure displays the optimal
frame bounds of $Z_{k}$ for each iteration step. Here we see exponential
convergence of the lower frame bound of $Z_{k}$ to zero.

\subsection{Comparison with other methods}

\begin{figure}
\includegraphics[%
  width=0.7\textwidth,
  keepaspectratio]{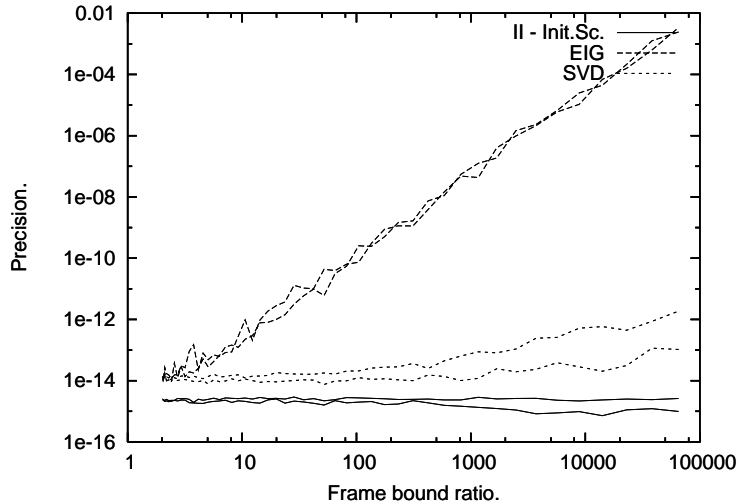}

\caption{\label{cap:numprec}The figure shows the numerical accuracy of three
different methods to compute the canonical tight window.The plot is
parametric in $w$. For each method, one line corresponds to a narrow
window, $w<1$, the other corresponds to a wide window, $w>1$. Almost
overlapping points on different lines corresponds to values $w_{1},w_{2}$
such that $w_{1}=1/w_{2}$. The error measure used is the dual lattice
norm. The Gabor frames used are the Gabor frames for $\mathbb{C}^{432}$,
$\left(\varphi_{w}^{D},18,18\right)$.\protect \\
}
\end{figure}

Figure \ref{cap:numprec} show a comparison of the numerical precision
of the algorithms for computing the canonical tight window compared
to the numerical precision of other standard methods. The stability
of the tight iterations proved in Subsection \ref{sub:Influence-of-out-of-space}
is clearly visible. The method based on computing eigenvalues deteriorates
quickly as the frame bound ratio increases, while the SVD behaves
much better. The eigenvalue method should not be used if the frame
bound ratio of the problem is unknown. An explanation for this is
that in the SVD method, the singular values are never considered,
they are simple set to $1$. Therefore, roundoff errors on the small
singular values do not affect the computation, in contrast to the
EIG method, where round-off errors on the smallest eigenvalues are
magnified because of the inversion of eigenvalues. For more details
on the stability of computing eigenvalues and singular values see
\cite{LAPACK99}.

The actual running time of the methods is determined by the flop count
for each method (see the previous section for details) and of how
fast the floating point operations can be executed by a computer.
We will not give timings of the iterative algorithms, because we have
not created optimal implementations of the algorithms, so timing them
makes little sense. We note, however, that the key ingredients in
the algorithms are FFTs of small length and matrix multiplications
of small size matrices. Fast implementations exists for both algorithms,
see \cite{whal00,FFTW05}. This makes it possible to create efficient
implementations of the iterative algorithms.

\subsection{Number of iterations}

\begin{figure}
\subfigure[Tight iterations.]{\includegraphics[%
  width=0.5\textwidth,
  keepaspectratio]{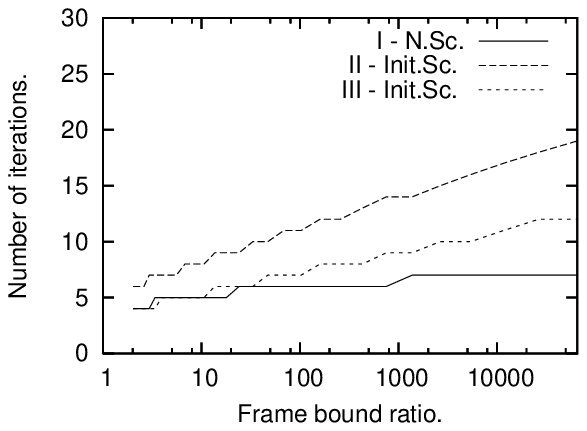}}\subfigure[Dual iterations.]{\includegraphics[%
  width=0.5\textwidth,
  keepaspectratio]{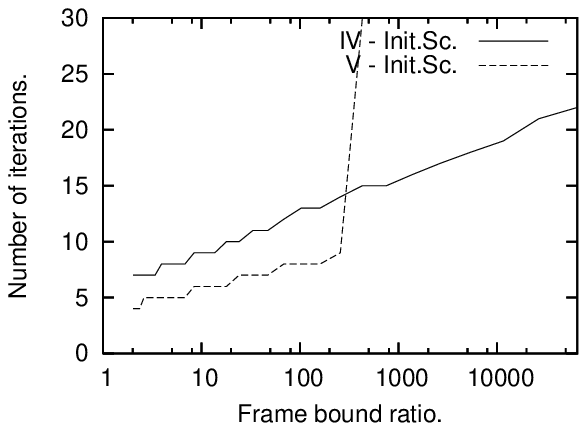}}

\caption{\label{cap:numits}The figure shows the number of iterations the
algorithms need in order to reach machine precision. The three algorithms
compared are I using norm scaling, and II-V using initial scaling
by (\ref{eq:Bbound}). The Gabor frames used are the Gabor frames
for $\mathbb{C}^{432}$, $\left(\varphi_{w}^{D},18,18\right)$.\protect \\
}
\end{figure}

To study how the number of necessary iterations depends on the frame
bound ratio of the initial Gabor frame, we have plotted the number
of iterations for the algorithms to converge, as a function of the
frame bound ratio of the Gabor frame. Figure \ref{cap:numits} shows
such a plot, using Gaussians to generate Gabor frames with varying
frame bound ratios. The jumps in the curves occur when, according
to the stopping criteria, an additional iteration step is necessary.
Even though algorithm III has cubic convergence, it is almost never
able to compete with I. The jump in the graph for V is due to the
magnification of round-off errors dominating the convergence, and
causing divergence. The same happens for algorithm IV, but for considerable
worse frame bound ratios (not visible on the graph). The graphs for
the hyperbolic secant look similar.

The number of necessary iterations might also depend on the size of
the matrix blocks appearing in the factorization. This issue is slightly
problematic to address, since creating a test problem involving bigger
matrices also means altering the frame bound ratio. To minimize this
effect, we have considered $p,q$ running through the Fibonacci numbers,
$2,3,5,8,\ldots$, such that $p/q\rightarrow\left(\sqrt{5}-1\right)/2$.
This creates a series of irreducible fractions $p/q$ while keeping
$p/q$ close to a certain number away from $1$. The result of the
test is that the number of necessary iterations seems to be completely
independent of the size of the matrix blocks! We have omitted the
graphs, as they are simply horizontal lines. For algorithm I, it is
proved in \cite{kenlaub92} that this is indeed the case.

\subsection{Choosing an initial scaling}

\begin{figure}
\subfigure[tight iterations]{\includegraphics[%
  width=0.5\textwidth,
  keepaspectratio]{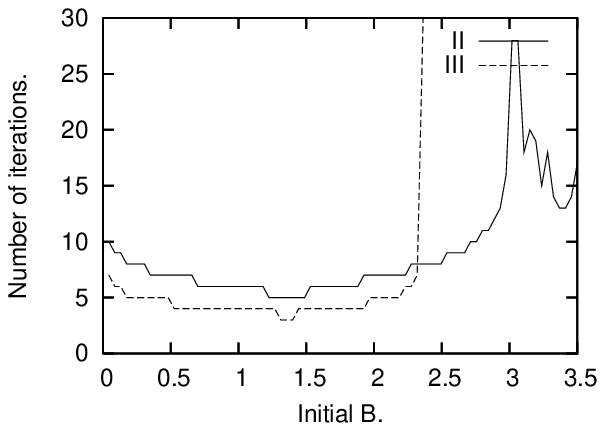}}\subfigure[dual iterations]{\includegraphics[%
  width=0.5\textwidth,
  keepaspectratio]{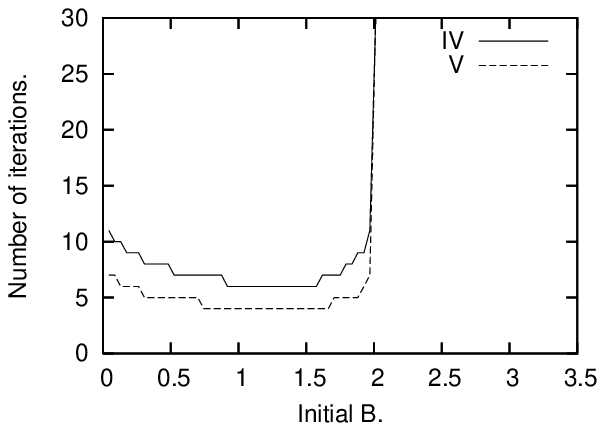}}

\caption{\label{cap:scaling}The figure shows the number of necessary iterations
to find the tight/dual window of a Gabor frame for $\mathbb{C}^{432}$,
$\left(\varphi_{1}^{D},18,18\right)$, as a function of the best upper
frame bounds of the initial window.\protect \\
}
\end{figure}

Figure \ref{cap:phi} shows that for each iteration type there is
a range of values of the upper frame bound of the scaled window, $B_{scaled}$,
that will guarantee convergence. Figure \ref{cap:scaling} shows an
example of the effect of prescaling the input window to obtain specific
values of $B_{scaled}$. As shown on Fig. \ref{cap:phi}, algorithm
III diverges if $B_{scaled}$ is larger than $7/3$. The dual iterations
IV and V diverge if $B_{scaled}>2$ and II has a chaotic behaviour
if $3<B_{scaled}<5$ and diverges if $B_{scaled}>5$ (not shown on
the plot).

The choice of $\hat{B}$ that minimizes the number of iterations is
to choose $\hat{B}$ according to Table \ref{cap:optscal}. An estimate
for this is difficult to calculate, as it involves an estimate for
the lower frame bound. Fortunately, as can be seen on Fig. \ref{cap:scaling}
there is a large region around the optimal scaling point, where only
1 or 2 extra iterations are needed.

\appendix

\section{\label{sec:Condition-A}A result on Condition A'}

\begin{prop}
Assume that $\left(g,a,b=1/M\right)$ is a Gabor frame that satisfies
condition A, where $a,M\in\mathbb{N}$. Also assume that $\varphi$
is analytic around $\left[A,B\right]$ and positive on $\left[A,B\right]$,
where $A>0$, $B<\infty$ are lower, upper frame bounds of $\left(g,a,b\right)$.
Finally, let $\gamma=\varphi\left(S\right)g$ where $S$ is the frame
operator corresponding to $\left(g,a,b\right)$. Then $g,\gamma$
satisfies the Condition A', i.e.,\begin{equation}
\sum_{j,l}\left|\left(g,\gamma_{j/b,l/a}\right)\right|<\infty.\end{equation}

\end{prop}
\begin{pf}
We have for $f,h\in L^{2}\left(\mathbb{R}\right)$ that \begin{eqnarray}
\sum_{n,m}\left(f,\gamma_{na,mb}\right)\left(g_{na,mb},h\right) & = & \sum_{n,m}\left(f,\left(\varphi\left(S\right)g\right)_{na.mb}\right)\left(g_{na,mb},h\right)\nonumber \\
 & = & \sum_{n,m}\left(\varphi\left(S\right)f,g_{na,mb}\right)\left(g_{na,mb},h\right)=\left(S\varphi\left(S\right)f,h\right).\end{eqnarray}

We know that $\left(\gamma,a,b\right)$is a Gabor frame. Now take
$f,h\in L^{2}\left(\mathbb{R}\right)$ such that $\left(f,a,b\right)$
and $\left(h,a,b\right)$ have finite upper frame bounds. Then by
the fundamental identity of Gabor analysis, see \cite[Subsecs. 1.4.1 and 1.4.2]{jans98},\begin{eqnarray}
\sum_{n,m}\left(f,\gamma_{na,mb}\right)\left(g_{na,mb},h\right) & = & \frac{1}{ab}\sum_{j,l}\left(g,\gamma_{j/b,l/a}\right)\left(f_{j/b,l/a},h\right),\label{eq:app-3}\end{eqnarray}
with absolute convergence on either side of (\ref{eq:app-3}) . Thus
$S\varphi\left(S\right)$ has the dual lattice representation\begin{eqnarray}
S\varphi\left(S\right) & = & \frac{1}{ab}\sum_{j,l}\left(g,\gamma_{j/lb,l/a}\right)U_{j,l}.\end{eqnarray}
Now let $\psi\left(s\right)=s\varphi\left(s\right)$. This $\psi$
is analytic around $\left[A,B\right]$ and positive on $\left[A,B\right]$.
By functional calculus of frame operators in the time-frequency domain,
see \cite[Sec. 8.3]{chrstr03}, there holds that $S\varphi\left(S\right)$
has also the dual lattice representation\begin{eqnarray}
S\varphi\left(S\right) & = & \sum_{j,l}\left(\psi\left(\frac{1}{ab}HH^{*}\right)\right)_{0,0;j,l}U_{j,l},\end{eqnarray}
where $H$ is the analysis operator with respect to the dual lattice,
defined for $f\in L^{2}\left(\mathbb{R}\right)$, by\begin{eqnarray}
Hf & = & \left(\left(f,g_{j/b,l/a}\right)\right)_{j,l\in\mathbb{Z}}.\end{eqnarray}
 It follows from the proof of \cite[Thm. 4.3]{chrstr03}, in particular
from uniform boundedness of (8.4.14) (with $\psi$ instead of $\varphi$),
that \begin{eqnarray}
\sum_{j,l}\left|\left(\psi\left(\frac{1}{ab}HH^{*}\right)\right)_{0,0;j,l}\right| & < & \infty.\end{eqnarray}
 By uniqueness of the coefficients in the dual lattice representation
(just consider a well-behaved $h$ such that $\left(h_{na,mb}\right)_{n,m\in\mathbb{Z}}$
is a tight frame, i.e. such that $U_{j,l}h$, $j,l\in\mathbb{Z}$,
is an orthogonal set of functions), it follows that $\sum_{j,l}\left|\left(g,\gamma_{j/b,l/a}\right)\right|<\infty$,
as required.
\end{pf}
\begin{note}
It is implicit in the statement and proof of \cite[Thm 4.3]{chrstr03}
that the $\gamma$ of the above result is such that condition A is
satisfied by $\left(\gamma,a,b\right)$.
\end{note}
\begin{ack}
The cooperation leading to this paper began in Vienna, 2005, during
the {}``Special Semester on Modern Methods of Time-Frequency Analysis'',
organized by H.G. Feichtinger and K. Gröchenig under the auspices
of the Erwin Schrödinger Institute for Mathematical Physics (ESI).
The authors express their thanks to the organizers and ESI.
\end{ack}
\bibliographystyle{abbrv}
\bibliography{ref}

\end{document}